\documentclass{article}
\usepackage{amsmath}
\usepackage{amssymb}
\usepackage{amsfonts}

\newtheorem{theorem}{Theorem}[section]

\newtheorem{corollary}[theorem]{Corollary}

\newtheorem{definition}[theorem]{Definition}

\newtheorem{lemma}[theorem]{Lemma}

\newtheorem{proposition}[theorem]{Proposition}
\newtheorem{remark}[theorem]{Remark}

\newenvironment{proof}[1][Proof]{\textbf{#1.}\;}{\vspace{-0,25 in} \begin{flushright}\rule{0.3em}{0.3em}
\end{flushright}}

\begin{document}

\title{A Functional Calculus \\[0pt]
For Quotient Bounded Operators \thanks{%
This work is supported by the M.Ed.C. grant C.N.B.S.S./2006 }}
\author{Sorin Mirel Stoian \\
%EndAName
University of Petro\c{s}ani }
\date{}
\maketitle\begin{abstract}

If $(X,\mathcal{P})$ is a sequentially complete locally convex space, then a quotient
bounded operator $T\in Q_{\mathcal{P}}(X)$ is regular (in the sense of
Waelbroeck) if and only if it is a bounded element (in the sense of Allan)
of algebra $Q_{\mathcal{P}}(X)$. The classic functional calculus for the bounded
operators on Banach space is generalized for bounded elements of algebra $Q_{%
\mathcal{P}}(X)$.
\end{abstract}


\section{ Introduction}

It is well-known that if $X$ is a Banach space and $L(X)$ is Banach algebra
of bounded operators on $X$, then formula

\begin{center}
$f(T)=\frac{1}{2\pi i}\int\limits_{\Gamma }f\left( z\right) R\left(
z,T\right) dz$,
\end{center}
( where $f$ is an analytic function on some neighborhood of $\sigma (T)$, $%
\Gamma $ is a closed rectifiable Jordan curve whose interior domain $D$ is
such that $\sigma (T)\subset D$, and $f$ is analytic on $D$ and continuous
on $D\cup \Gamma $) defines a homomorphism $f\rightarrow f(T)$ from the set
of all analytic functions on some neighborhood of $\sigma (T)$ into $L(X) $,
with very useful properties.

Through this paper all locally convex spaces will be assumed Hausdorff, over
complex field $\mathbb{C}$, and all operators will be linear. If $X$ and $Y$
are topological vector spaces we denote by $L(X,Y)$ ($\mathcal{L}(X,Y)$) the
algebra of linear operators (continuous operators) from $X$ to $Y$.

Any family $\mathcal{P}$ of seminorms which generate the topology of locally
convex space $X$ (in the sense that the topology of $X$ is the coarsest with
respect to which all seminorms of $\mathcal{P}$ are continuous) will be
called a calibration on $X$.The set of all calibrations for $X$ is denoted
by $\mathcal{C}(X)$ and the set of all principal calibration by $\mathcal{C}%
_{0}(X)$.

If $(X,\mathcal{P})$ is a locally convex algebra and each seminorms $p\in 
\mathcal{P}$ \ is submultiplicative then we said that $(X,\mathcal{P})$ is
l.m.c.-algebra.

On a family of seminorms on a linear space X we define the relation ,, $\leq 
$' of partially ' ordered by

\begin{center}
$p\leq q\Leftrightarrow p\left( x\right) \leq q\left( x\right) $, $\left(
\forall \right) x\in X$.
\end{center}

A family of seminorms is preordered by relation ,, $\prec $'', where

\begin{center}
$p\prec q\Leftrightarrow $ there exists some $r>0$ such that $p\left(
x\right) \leq rq\left( x\right) $, for all $x\in X$.
\end{center}

If $p\prec q$ and $q\prec p$, we write $p\approx q$.

\begin{definition}
Two families $\mathcal{P}_{1}$ and $\mathcal{P}_{2}$ of seminorms on a
linear space are called $Q $-equivalent ( denoted $\mathcal{P}_{1}\approx 
\mathcal{P}_{2}$) provided:
\begin{enumerate}
\item 
for each $\ p_{1}\in $ $P_{1}$ there exists $p_{2}\in \mathcal{P}_{2}$
such that $p_{1}\approx p_{2}$;
\item 
for each $p_{2}\in P_{2}$ there exists $p_{1}\in \mathcal{P}_{1}$ such
that $p_{2}\approx p_{1}$.
\end{enumerate}
\end{definition}

It is obvious that two $Q$-equivalent and separating families of seminorms
on a linear space generate the same locally convex topology.

Similar to the norm of an operator on a normed space we define the mixed
operator seminorm of an operator between locally convex spaces. If $(X,%
\mathcal{P})$, $(Y,\mathcal{Q})$ are locally convex spaces, then for each $%
p\in \mathcal{P}$ and $q\in \mathcal{Q}$ the application $m_{pq}:L(X,Y)\rightarrow \mathbb{R} \cup \left\{
\infty \right\} $, defined by
\begin{equation*}
m_{pq}(T)=\sup\limits_{p\left( x\right) \neq 0}\frac{q\left( Tx\right) }{%
p\left( x\right) },\left( \forall \right) T\in L(X,Y).
\end{equation*}
is called the mixed operator seminorm of $T$ associated with $p$ and $q$.
When $X=Y$ and $p=q$ we use notation $\hat{p}=m_{pp}$.

\begin{lemma}[\cite{tr}]
{\label{lemma:op.semin.}}If $(X,\mathcal{P}),(Y,\mathcal{Q})$ are locally
convex spaces and $T\in L(X,Y)$, then
\begin{enumerate}
\item 
$m_{pq}(T)=\sup\limits_{p\left( x\right) =1}q\left( Tx\right)
=\sup\limits_{p\left( x\right) \leq 1}q\left( Tx\right) ,\left( \forall
\right) p\in P,\left( \forall \right) q\in Q;$
\item 
$\ q\left( Tx\right) \leq m_{pq}\left( T\right) p\left( x\right) $, $%
\left( \forall \right) x\in X$, whenever $m_{pq}(T)$ $<\infty $.
\item
$m_{pq}(T)=\inf \left\{ M>0\left\vert q\left( Tx\right) \leq Mp\left(
x\right) \right. ,\left( \forall \right) x\in X\right\} $, whenever \\ $%
m_{pq}(T)<\infty $.
\end{enumerate}
\end{lemma}

\begin{definition}
Let X be a locally convex space. An operator $T\in L(X)$ is:
\begin{enumerate}
\item 
quotient bounded operator with respect to a calibration $\mathcal{P}%
\in \mathcal{C}(X)$ if for every seminorm $p\in \mathcal{P}$ there exists $%
c_{p}>0$ such that 
\begin{equation*}
p\left( Tx\right) \leq c_{p}p\left( x\right) ,\left( \forall \right) x\in X.
\end{equation*}
\item
universally bounded with respect to the calibration $\mathcal{P}\in 
\mathcal{C}(X)$ if there exists $c_{0}>0$ such that 
\begin{equation*}
p(Tx)\leq c_{0}p\left( x\right) ,\left( \forall \right) x\in X,\left(
\forall \right) p\in \mathcal{P}.
\end{equation*}
\item 
locally bounded if maps some zero neighborhood in a bounded set, i.e.
there exists some seminorms $p\in \mathcal{P}$ such that we have $%
m_{pq}(T)<\infty $, for each $q\in \mathcal{P}$ .
\end{enumerate}
\end{definition}

\begin{remark}
It is obvious that $B_{\mathcal{P}}(X)\subset Q_{\mathcal{P}}(X)\subset 
\mathcal{L}(X)$.
\end{remark}

If X is a locally convex space and $\mathcal{P}$ $\in C(X)$, then for every $%
p\in \mathcal{P}$ the application $\hat{p}:Q_{\mathcal{P}}(X)\rightarrow 
\mathbf{R}$ defined by 
\begin{equation*}
\hat{p}(T)=\inf \{~r>0~{|~}p(Tx)\leq r~p\left( x\right) ,\left( \forall
\right) x\in X~\},
\end{equation*}%
is a submultiplicative seminorm on $Q_{\mathcal{P}}(X)$, satisfying $\hat{p}%
(I)=1$. We denote by $\hat{\mathcal{P}}$ the family $\{~\hat{\mathcal{P}}~{|~%
}p\in \mathcal{P}\}$.

\begin{lemma}[\cite{st}]
If $X$ is a sequentially complete convex space, then $Q_{\mathcal{P}}(X)$ is
a sequentially complete m-convex algebra for all $\mathcal{P}\in \mathcal{C}%
(X) $.
\end{lemma}

Let $(X,\mathcal{P})$ be a locally convex space and $T\in Q_{\mathcal{P}}(X)$. We say that $T$ is a bounded element of the algebra $Q_{\mathcal{P}}(X)$ if
it is bounded element in the sense of G.R.Allan $\left[ 1\right] $, i.e some
scalar multiple of it generates a bounded semigroup. The class of the bounded
elements of $Q_{\mathcal{P}}(X)$ is denoted by $(Q_{\mathcal{P}}(X))_{0}$.

\begin{definition}
Let $(X,\mathcal{P})$ be a locally convex space.
\begin{enumerate}
\item 
If $T\in Q_{\mathcal{P}}(X)$ ($T\in (Q_{\mathcal{P}}(X))_{0}$, respectively $T\in B_{\mathcal{P}}(X)$) we said that $%
\alpha \in \mathbb{C}$ is in the resolvent set $\rho (Q_{\mathcal{P}},T)$ ($%
\rho (Q_{\mathcal{P}}^{0},T)$, respectively $\rho (B_{\mathcal{P}},T)$) if
there exists $(\alpha I-T)^{-1}\in Q_{\mathcal{P}}(X)$ ($(\alpha
I-T)^{-1}\in (Q_{\mathcal{P}}(X))_{0}$, respectively $(\alpha I-T)^{-1}\in
B_{\mathcal{P}}(X)$). The spectral set $\sigma (Q_{\mathcal{P}},T)$ ($\sigma
(Q_{\mathcal{P}}^{0},T)$, respectively $\sigma (B_{\mathcal{P}},T)$) will be
the complement set of $\rho (Q_{\mathcal{P}},T)$ ($\rho (Q_{\mathcal{P}%
}^{0},T)$, respectively $\rho (B_{\mathcal{P}},T)$).
\item 
Let $T$ be a locally bounded operator on a locally convex space $X$. We
say that $\lambda \in \rho _{lb}(T)$ if there exists a scalar $\alpha $ and
a locally bounded operator $S$ on $X$ such that $(\lambda I-T)^{-1}=\alpha
I+S$. The spectral set $\sigma _{lb}(T$) is defined to be
the complement of the resolvent set $\rho _{lb}(T$).
\end{enumerate}
\end{definition}

\begin{remark}
\begin{enumerate}
\item 
It is obvious that we have the following inclusions 
\begin{equation*}
\sigma (T)\subset \sigma (Q_{\mathcal{P}},T)\subset \sigma (B_{\mathcal{P}%
},T).
\end{equation*}
\item 
The set $\rho _{lb}(T)$ will be the spectrum of $T$ in the algebra $%
\mathcal{LB}_{0}(X)$. 
\end{enumerate}
\end{remark}

\begin{definition}
Let $X$ be a locally convex space. A sequence $(S_{n})_{n}\subset \mathcal{LB%
}(X)$ converges uniformly to zero on some zero-neighborhood if for each
principal calibration $\mathcal{P}\in \mathcal{C}_{0}(X)$ there exists some
seminorm $p\in \mathcal{P}$ such that for every $q\in \mathcal{P}$ and every 
$\epsilon >0$ there exists an index $n_{q,\epsilon }\in \mathbf{N}$, with
the property 
\begin{equation*}
m_{pq}(S_{n})<\epsilon ,\left( \forall \right) n\geq n_{q,\epsilon }.
\end{equation*}%
A family $G\subset \mathcal{LB}(X)$ is uniformly bounded on some
zero-neighborhood if there exists some seminorm $p\in \mathcal{P}$ such that
for every $q\in \mathcal{P}$ there exists $\epsilon _{q}>0$ with the
property 
\begin{equation*}
m_{pq}(S)<\epsilon _{q},\left( \forall \right) S\in G.
\end{equation*}
\end{definition}

\begin{proposition}[\cite{jo}]
{\label{proposition:jo1}}Let $X$ be a locally convex space and $\mathcal{P}%
\in \mathcal{C}(X)$.
\begin{enumerate}
\item 
$Q_{\mathcal{P}}(X)$ is a unitary subalgebra of $\mathcal{L}(X)$;
\item 
$Q_{\mathcal{P}}(X)$ is a unitary l.m.c.-algebra with respect to the
topology determined by $\hat{\mathcal{P}}$;
\item 
If $\mathcal{P}^{\prime }\in \mathcal{C}(X)$ such that $\mathcal{P}%
\approx \mathcal{P}^{\prime }$, then $Q_{\mathcal{P}^{\prime }}(X)=Q_{%
\mathcal{P}}(X)$; moreover the $\hat{\mathcal{P}}$-topology and $\hat{%
\mathcal{P}}^{\prime }$-topology coincide;
\item 
The topology generated by $\hat{\mathcal{P}}$ on $Q_{\mathcal{P}}(X)$ is finer than the topology of uniform convergence on bounded subsets of $X$.
\end{enumerate}
\end{proposition}

\begin{definition}
If $(X,\mathcal{P})$ is a locally convex space and $T\in Q_{\mathcal{P}}(X)$
we denote by $r_{\mathcal{P}}(T)$ the radius of boundness of operator $T$ in 
$Q_{\mathcal{P}}(X)$, i.e. 
\begin{equation*}
r_{\mathcal{P}}(T)=\inf \{\alpha >0\;\,\mid \;\,\alpha ^{-1}T\text{
generates a bounded semigroup in }Q_{\mathcal{P}}(X)\}.
\end{equation*}
\end{definition}

We have said that $r_{\mathcal{P}}(T)$ is the $\mathcal{P}$-spectral radius
of the operator $T$.

\begin{proposition}[\cite{al}]
{\label{proposition:al}}If $X$ is a locally convex space and $\mathcal{P}\in 
\mathcal{C}(X)$, then for each $T\in Q_{\mathcal{P}}(X)$ we have 
\begin{equation*}
r_{\mathcal{P}}(T)=\sup \{\;\,\underset{n\rightarrow \infty }{\limsup }%
\left( \hat{p}\left( T^{n}\right) \right) ^{1/n}\mid \;p\in \mathcal{P}\}.
\end{equation*}
\end{proposition}

\begin{proposition}[\cite{st}]{\label{proposition:st1}} 
If $X$ is a locally convex space and $\mathcal{P}\in \mathcal{C}(X)$, then
for each $T\in Q_{\mathcal{P}}(X)$ we have:
\begin{enumerate}
\item 
$r_{\mathcal{P}}(T)\geq 0$ and $r_{\mathcal{P}}(\lambda T)=\mid
\lambda \mid r_{\mathcal{P}}(T),\left( \forall \right) \lambda \in \mathbb{C}
$, where by convention $0\infty =\infty $;
\item 
$r_{\mathcal{P}}(T)<+\infty $ if and only if $T\in (Q_{\mathcal{P}%
}(X))_{0}$;
\item 
$r_{\mathcal{P}}(T)=\inf \left\{ \lambda >0\mid
\;\lim\limits_{n\rightarrow \infty }\frac{T^{n}}{\lambda ^{n}}=0\;\right\} $;
\item 
$r_{\mathcal{P}}(T)=\sup \{\;\,\lim\limits_{n\rightarrow \infty }\left( \hat{p}\left(
T^{n}\right) \right) ^{1/n}\mid \;p\in \mathcal{P}\}=\newline
=\sup\{\;\,\inf\limits_{n\geq 1}\left( \hat{p}\left( T^{n}\right) \right)
^{1/n}\mid \;p\in \mathcal{P}\};$
\end{enumerate}
\end{proposition}

\begin{proposition}[\cite{st}]
Let $X$ be a locally convex space and $\mathcal{P}$ $\in \mathcal{C}(X)$.
\begin{enumerate}
\item 
If $T\in (Q_{\mathcal{P}}(X))_{0}$, then 
\begin{equation*}
\lim\limits_{n\rightarrow \infty }\frac{T^{n}}{\lambda ^{n}}=0, \left(
\forall \right) \mid \lambda \mid >r_{\mathcal{P}}(T);
\end{equation*}
\item 
If $T\in (Q_{\mathcal{P}}(X))_{0}$ and $0<\mid \lambda \mid <r$ $_{%
\mathcal{P}}(T)$, then the set $\left\{ \frac{T^{n}}{\lambda ^{n}}\right\}
_{n\geq 1}$ is unbounded.
\item 
For each $T\in Q_{\mathcal{P}}(X)$ and every $n$ $>0$ we have $r_{%
\mathcal{P}}(T^{n})=r_{\mathcal{P}}(T)^{n}$.
\end{enumerate}
\end{proposition}

\begin{proposition}[\cite{st}]
If $(X,\mathcal{P})$ is a sequentially complete locally convex space, then $Q_{%
\mathcal{P}}(X)$ is a sequentially complete l.m.c.-algebra for all $\mathcal{%
P}\in \mathcal{C}(X)$.
\end{proposition}

\begin{proposition}[\cite{st}]
{\label{proposition:newman}} Let $X$ be a sequentially complete locally
convex space and $\mathcal{P}\in \mathcal{C}(X)$. If $T\in (Q_{\mathcal{P}%
}(X))_{0}$ and $\mid \lambda \mid >r_{\mathcal{P}}(T)$, then the Neumann
series $\sum\limits_{n=0}^{\infty }\frac{{\ T}^{n}}{\lambda ^{n+1}}$
converges to $R\left( \lambda ,T\right) $(in $Q_{\mathcal{P}}(X)$) and $%
R\left( \lambda ,T\right) \in Q_{\mathcal{P}}(X)$.
\end{proposition}

\begin{proposition}[\cite{st}]{\label{proposition:raza spectrala}}
Let X be a sequentially complete locally
convex space and $\mathcal{P}\in C(X)$. If $T\in Q_{\mathcal{P}}(X)$, then $%
\left\vert \sigma (Q_{\mathcal{P}},T)\right\vert =r_{\mathcal{P}}(T).$
\end{proposition}

\begin{definition}
If $(X,\mathcal{P})$ is a locally convex space and $T\in (Q_{\mathcal{P}%
}(X))_{0}$ we denote by $r_{\mathcal{P}}^{0}(T)$ the radius of boundness of the
operator $T$ in $(Q_{\mathcal{P}}(X))_{0}.$
\end{definition}

We say that $r_{\mathcal{P}}^{0}(T)$ is the $\mathcal{P}$-spectral radius of
the operator $T$ in algebra $(Q_{\mathcal{P}}(X))_{0}$.

\begin{remark}
If $(X,\mathcal{P})$ is a locally convex space and $T\in (Q_{\mathcal{P}}(X))_{0}$
then $r_{\mathcal{P}}^{0}(T)=r_{\mathcal{P}}(T)$. Moreover, $r_{\mathcal{P}%
}^{0}(T)$ has all the properties of the spectral radius $r_{\mathcal{P}}(T)$
presented above.
\end{remark}

\begin{lemma}[\cite{st}]
If $\mathcal{P}$ a calibration on $X$, then $B_{\mathcal{P}}(X)$ is a
unitary normed algebra with respect to the norm $\left\Vert {\bullet }%
\right\Vert _{\mathcal{P}}$ defined by 
\begin{equation*}
\left\Vert T\right\Vert _{\mathcal{P}}=\sup \{m_{pp}(T)\;|\;p\in \mathcal{P}\},\left(
\forall \right) T\in B_{\mathcal{P}}(X).
\end{equation*}
\end{lemma}

\begin{corollary}[\cite{st}]
If $\mathcal{P}\in\mathcal{C}(X)$, then for each $T\in B_\mathcal{P}(X)$ we
have 
\begin{equation*}
\left\Vert T\right\Vert _\mathcal{P}=\inf \{M>0~|~~p\left( Tx\right) \leq
Mp\left( x\right) ,\left( \forall \right) x\in X,\left( \forall \right)
~p\in \mathcal{P}\}.
\end{equation*}
\end{corollary}

\begin{proposition}[\cite{jo}]
Let $X$ be a locally convex space and $\mathcal{P}\in \mathcal{C}(X)$. Then:

\begin{enumerate}
\item $B_{\mathcal{P} }(X)$ is a subalgebra of $\mathcal{L}(X)$;

\item ($B_{\mathcal{P} }(X), \left\Vert {\bullet}\right\Vert _{\mathcal{P}}$%
) is unitary normed algebra;

\item for each $\mathcal{P}^{\prime }\in \mathcal{C}(X)$ with the property $%
\mathcal{P}\sim \mathcal{P}^{\prime }$, we have 
\begin{equation*}
B_{\mathcal{P}}(X)=B_{\mathcal{P}^{\prime }}(X)\;\text{and}\;\left\Vert {%
\bullet }\right\Vert _{\mathcal{P}}=\left\Vert {\bullet }\right\Vert _{%
\mathcal{P}^{\prime }}.
\end{equation*}
\end{enumerate}
\end{proposition}

\begin{proposition}[\cite{ch}]
Let $X$ be a locally convex space and $\mathcal{P}\in \mathcal{C}(X)$. Then:
\begin{enumerate}
\item the topology given by the norm $\left\Vert {\bullet}\right\Vert $ $_{%
\mathcal{P}}$ on the algebra $B_{\mathcal{P} }(X)$ is finer than the
topology of uniform convergence;
\item 
if $(T_{n})_{n}$ is a Cauchy sequences in $(B_{\mathcal{P}
}(X),\left\Vert {\bullet}\right\Vert _{\mathcal{P} })$ which converges punctually to an
operator $T$, we have $T\in B_{\mathcal{P} }(X)$;
\item 
the algebra $(B_{\mathcal{P} }(X),\left\Vert {\bullet}\right\Vert _{%
\mathcal{P}})$ is complete if $X$ is sequentially complete.
\end{enumerate}
\end{proposition}

\begin{proposition}[\cite{jo}]{\label{proposition:jo2}} 
Let $(X,\mathcal{P})$ be a locally convex space. An operator $T\in Q_{\mathcal{P}}(X)$ is bounded in the algebra $Q_{\mathcal{%
P}}(X)$ if and only if there exists some calibration $\mathcal{P}^{\prime
}\in \mathcal{C}(X)$ such that $\mathcal{P}\approx \mathcal{P}^{\prime }$
and $T\in B_{\mathcal{P}^{\prime }}(X)$.
\end{proposition}

\begin{proposition}[\cite{ch}]
If $(X,\mathcal{P})$ is a sequentially complete locally convex space and $T\in B_ \mathcal{P}(X)$,
then the set $\sigma (B_{\mathcal{P}},T)$ is compact.
\end{proposition}

\begin{lemma}{\label{lemma:1}} 
If $\mathcal{P}$ is a calibration on a locally convex
space $X$, then 
\begin{equation*}
\left\Vert T\right\Vert _{\mathcal{P}}=\sup \{\hat{p}(T)\left\vert \;p\in
\right. \mathcal{P}\},\left( \forall \right) T\in B_{\mathcal{P}}(X).
\end{equation*}
\end{lemma}

\begin{corollary}
If $X$ is a locally convex algebra $\mathcal{P}\in \mathcal{C}(X)$, then for
each $T\in B_{\mathcal{P}}(X)$ the inequality $r_{\mathcal{P}}(T)\leq
\left\Vert T\right\Vert _{\mathcal{P}^{\prime }}$ holds for each calibration 
$\mathcal{P}^{\prime }\in \mathcal{C}(X)$ such that $\mathcal{P}\approx 
\mathcal{P}^{\prime }$ and $T\in B_{\mathcal{P}^{\prime }}(X)$.
\end{corollary}

\begin{lemma}[\cite{st1}]{\label{proposition:univ.bound.1}} 
Let $X$ be a sequentially complete
locally convex space and $\mathcal{P}\in \mathcal{C}(X)$. Then, for each $%
T\in B_{\mathcal{P}}(X)$ we have 
\begin{equation*}
\mid \sigma (B_{\mathcal{P} },T)\mid \leq \underset{n\rightarrow \infty }{%
\liminf }\left\Vert T^{n}\right\Vert _{\mathcal{P} }^{1/n}\leq r(B_{\mathcal{%
P} },T)
\end{equation*}
\end{lemma}

\begin{proposition}[\cite{gi}]{\label{proposition:gi2}} 
Let $(X,\mathcal{P}=\left( p_{\alpha }\right) _{\alpha \in \Lambda })$. If $%
T\in Q_{\mathcal{P}}(X)$, such that the set $\sigma (Q_{\mathcal{P}},T)$ is
bounded, then there exists some calibration $\mathcal{P}^{\prime }=\left(
p_{\alpha }^{\prime }\right) _{\alpha \in \Lambda }\in \mathcal{C}(X)$ with
the property:

(A) For all $\alpha \in \Lambda $ there exists $m_{\alpha },M_{\alpha }>0$
such that 
\begin{equation*}
m_{\alpha }p_{\alpha }\left( x\right) \leq p_{\alpha }^{^{\prime }}\left(
x\right) \leq M_{\alpha }p_{\alpha }\left( x\right) ,\left( \forall \right)
x\in X.
\end{equation*}%
Moreover, $T\in B_{\mathcal{P}^{\prime }}(X)$.
\end{proposition}

\begin{definition}
Given a linear operator $T$ on a topological vector space X, we consider 
\begin{equation*}
r_{lb}(T)=\inf \left\{ \nu >0\left\vert \dfrac{{T}^{n}}{\nu ^{n}}\rightarrow
0\text{ uniformly on some zero neigborhood }\right. \right\}
\end{equation*}
\end{definition}

\begin{lemma}[\cite{k2}]{\label{lemma:k2}} 
If $T_{1}$ and $T_{2}$ are locally bounded operators on $%
X $, then there exists a calibration $\mathcal{P}^{\prime }$ on $X$ such
that $T_{1},T_{2}\in B_{\mathcal{P}^{\prime }}(X)$.
\end{lemma}

\begin{proposition}[\cite{st1}]{\label{proposition:spectru lb 1}} 
If $X$ is a sequentially complete locally
convex space and $T$ is locally bounded, then 
$\sigma (T)=\sigma _{lb}(T)$.
\end{proposition}

\begin{proposition}[\cite{st1}]
{\label{proposition:raza spectrala lb}} If $X$ is a sequentially complete
locally convex space and $T$ is locally bounded, then 
\begin{equation*}
r_{lb}(T)=\left\vert \sigma (T)\right\vert =\left\vert \sigma
_{lb}(T)\right\vert .
\end{equation*}
\end{proposition}

\begin{proposition}[\cite{st1}]
{\label{proposition:spectru lb}}If $X$ is a sequentially complete locally
convex space and $T$ is locally bounded, then $\sigma (T)$ is compact.
\end{proposition}

\begin{corollary}[\cite{st1}]
Let X be a sequentially complete locally convex space and $T\in \mathcal{LB}%
(X)$. If $\lambda \in \rho (T)$ and $d\left( \lambda \right) $ is the
distance from $\lambda $ to the set $\sigma (T)$, then 
\begin{equation*}
\left\Vert \left( \lambda I-T\right) ^{-1}\right\Vert _{\mathcal{P}}\geq 
\frac{1}{d\left( \lambda \right) },
\end{equation*}%
whenever $\mathcal{P}\in \mathcal{C}(X)$, such that $(\lambda I-T)^{-1},T\in
B_{\mathcal{P}}(X)$.
\end{corollary}

\begin{definition}{\label{corollary:1}}
{\label{definition:sw}} Let $(X,\mathcal{P})$ be a locally convex space. The
Waelbroeck resolvent set of an operator $T\in Q_{\mathcal{P}}(X)$, denoted
by $\rho _{W}(Q_{\mathcal{P}},T)$, is the subset of elements of $\lambda
_{0}\in \mathbb{C} _{\infty }=\mathbb{C}\cup \left\{ \infty \right\} $, for
which there exists a neighborhood $V\in \mathcal{V}_{(\lambda _{0})}$ such
that:
\begin{enumerate}
\item 
the operator $\lambda I-T$ is invertible in $Q_{\mathcal{P}}(X)$ for
all $\lambda \in V\backslash \{\infty \}$
\item 
the set $\{~\left( ~\lambda I-T~\right) ^{-1}|~\lambda \in V\backslash
\{\infty \}~\}$ is bounded in $Q_{\mathcal{P}}(X)$.
\end{enumerate}

The Waelbroeck spectrum of $T$, denoted by $\sigma _{W}(Q_{\mathcal{P}},T)$,
is the complement of the set $\rho _{W}(Q_{\mathcal{P}},T)$ in $\mathbb{C}%
_{\infty }$. It is obvious that $\sigma (Q_{\mathcal{P}},T)\subset \sigma
_{W}(Q_{\mathcal{P}},T)$.
\end{definition}

\begin{remark}
Let $(X,\mathcal{P})$ be a locally convex space.
If $T\in Q_{\mathcal{P}}(X))$, then  $\rho _{W}(Q_{\mathcal{P}^{\prime}},T)=\rho _{W}(Q_{\mathcal{P}},T)$ for all  $\mathcal{P}^{\prime}\in \mathcal{C}(X)$ such that $\mathcal{P}\approx \mathcal{P}^{\prime }$.
\end{remark}

\begin{definition}
{\label{definition:reg}} Let $(X,\mathcal{P})$ be a locally convex space. An
operator $T\in Q_{\mathcal{P}}(X)$ is regular if $\ \infty \notin\sigma
_{W}(Q_{\mathcal{P}},T)$, i.e. there exists some $t>0$ such that:

\begin{enumerate}
\item the operator $\lambda I-T$ is invertible in $Q_{\mathcal{P}}(X)$, for
all $\mid \lambda \mid >t$

\item the set $\left\{ R\left( \lambda ,T\right) \mid \mid \lambda \mid
>t\right\} $ is bounded in $Q_{\mathcal{P}}(X)$
\end{enumerate}
\end{definition}

\section{Bounded operators in $Q_{\mathcal{P}}(X)$}

In the next sections we assume that $X$ will be sequentially complete
locally convex space.

\begin{lemma}
{\label{lemma:qb0}} Let $(X,\mathcal{P})$ be a locally convex space and $%
T\in (Q_{\mathcal{P}}(X))_{0}$ such that $r_{\mathcal{P}}(T)<1$. Then the
operator $I-T$ is invertible and $I-T=\sum\limits_{n=0}^{\infty }T^{n}$.
\end{lemma}
\begin{proof}
Assume that $r_{\mathcal{P}}(T)<t<1$. From proposition \ref{proposition:al}
results that 
\begin{equation*}
\underset{n\rightarrow \infty }{\limsup }\left( \hat{p}\left( T^{n}\right)
\right) ^{1/n}<t,\left( \forall \right) \;p\in \mathcal{P},
\end{equation*}
so for each $\,\;p\in \mathcal{P}$ there exists $n_{p}\in \mathbb{N}$ such that
\begin{equation*}
(\hat{p}\left( T^{n}\right))^{1/n} \leq \underset{n\geq n_{p}}{\sup }\left( \hat{p}%
\left( T^{n}\right) \right) ^{1/n}<t,\left( \forall \right) \;n\geq n_{p}.
\end{equation*}

This relation implies that the series $\sum\limits_{n=0}^{\infty }\hat{p}%
\left( T^{n}\right) $ converges, so 
\begin{equation*}
\underset{n\rightarrow \infty }{\lim }\hat{p}\left( T^{n}\right) =0,\left(
\forall \right) \;p\in \mathcal{P},
\end{equation*}
therefore $\underset{n\rightarrow \infty }{\lim }T^{n}=0$. Since the algebra 
$Q_{\mathcal{P}}(X)$ is sequentially complete results that the series $%
\sum\limits_{n=0}^{\infty }T^{n}$ converges. Moreover, 
\begin{equation*}
(I-T)\sum\limits_{n=0}^{m}T^{n}=\sum\limits_{n=0}^{m}T^{n}(I-T)=I-T^{m+1},
\end{equation*}
so
\begin{equation*}
(I-T)\sum\limits_{n=0}^{\infty }T^{n}=\sum\limits_{n=0}^{\infty
}T^{n}(I-T)=I,
\end{equation*}
which implies that $I-T$ is invertible and $I-T=\sum\limits_{n=0}^{\infty
}T^{n}$.
\end{proof}

\begin{lemma}
{\label{lemma:qb1}} Let $(X,\mathcal{P})$ be a sequentially complete locally
convex space. If $T\in (Q_{\mathcal{P}}(X))_{0}$ then
\begin{enumerate}
\item 
the application $\lambda \rightarrow R(\lambda ,T)$ is holomorphic on $%
\rho _{W}(Q_{\mathcal{P}},T)$;
\item 
$\frac{d^{n}}{d\lambda ^{n}}R(\lambda ,T)=(-1)^{n}n!R(\lambda
,T)^{n+1} $, for every $n\in \mathbb{N}$;
\item 
$\lim\limits_{\vert \lambda \vert \rightarrow \infty }R(\lambda ,T)=0$ and $%
\lim\limits_{\vert \lambda \vert \rightarrow \infty }R(1,\lambda
^{-1}T)=\lim\limits_{\vert \lambda \vert \rightarrow \infty }\lambda R(1,T)=I$;
\item 
$\sigma _{W}(Q_{\mathcal{P}},T)\neq \varnothing .$
\end{enumerate}
\end{lemma}
\begin{proof}
1) If $\lambda _{0}\in \rho _{W}(Q_{\mathcal{P}},T)$ then there exists $V\in 
\mathcal{V}_{(\lambda _{0})}$ with the properties (1) and (2) from
definition (\ref{definition:sw}). Since for every $\lambda \in V\backslash
\{\infty \}$ we have 
\begin{equation*}
R(\lambda ,T)-R(\lambda _{0},T)=(\lambda _{0}-\lambda )R(\lambda
,T)R(\lambda _{0},T)
\end{equation*}%
and the set $\{R(\lambda ,T)|~\lambda \in V\backslash \{\infty \}\}$ is
bounded in $Q_{\mathcal{P}}(X)$ results that the application $\lambda
\rightarrow R(\lambda ,T)$ is continuous in $\lambda _{0}$, so 
\begin{equation*}
\lim\limits_{\lambda \rightarrow \lambda_{0} }\frac{R(\lambda ,T)-R(\lambda_{0}
,T)}{\lambda -\lambda _{0}}=-R^{2}(\lambda _{0},T)
\end{equation*}

If $\lambda _{0}=\infty $ then, there exists some neighborhood $V\in 
\mathcal{V}_{(\infty )}$ such that the application $\lambda \rightarrow
R(\lambda ,T)$ is defined and bounded on $V\backslash \{\infty \}$.
Moreover, this application it is holomorphic and bounded on $V\backslash
\{\infty \}$, which implies that it is holomorphic at $\infty $.

Therefore, the application $\lambda \rightarrow R(\lambda ,T)$ is
holomorphic on $\rho _{W}(Q_{\mathcal{P}},T)$.

2) Results from the proof of (1).

3) For each $\lambda \in \rho _{W}(Q_{\mathcal{P}},T)$ we have 
\begin{equation*}
\lambda ^{-1}(I+TR(\lambda ,T))(\lambda I-T)=I,
\end{equation*}%
so 
\begin{equation}
{\label{equation:21}}R(\lambda ,T)=\lambda ^{-1}(I+TR(\lambda ,T)).
\end{equation}

If $V\in \mathcal{V}_{(\lambda _{0})}$ satisfies the condition of the
definition (\ref{definition:sw}), then the set 
\begin{equation*}
\{TR(\lambda ,T)|~\lambda \in V\backslash \{\infty \}\}
\end{equation*}%
is bounded, so from relation (\ref{equation:21}) results that $%
\lim\limits_{\vert\lambda\vert \rightarrow \infty }R(\lambda ,T)=0.$

From equality $R(\lambda ,T)=\lambda ^{-1}R(1,\lambda ^{-1}T), \lambda\neq 0$, and relation (%
\ref{equation:21}) results that 
\begin{equation*}
R(1,\lambda ^{-1}T)=I+TR(\lambda ,T),
\end{equation*}%
so 
\begin{equation*}
\lim\limits_{\vert \lambda\vert  \rightarrow \infty }R(1,\lambda
^{-1}T)=\lim\limits_{\vert \lambda \vert \rightarrow \infty }(I+TR(\lambda ,T))=I
\end{equation*}%
4) Assume that $\sigma _{W}(Q_{\mathcal{P}},T)=\varnothing $. Then the
application $\lambda \rightarrow R(\lambda ,T)$ is holomorphic on $\mathbb{C}
$ and converges to 0 at infinity. From Liouville Theorem results that $%
R(\lambda ,T)=0,~\left( \forall \right) \lambda \in \mathbb{C}$. Therefore, $%
I=(\lambda I-T)R(\lambda ,T)=0$, which is not true.
\end{proof}

\begin{proposition}
{\label{proposition:qb1}} Let $(X,\mathcal{P})$ be a locally convex space.
Then an operator $T$ $\in Q_{\mathcal{P}}(X)$ is regular if and only if $%
T\in (Q_{\mathcal{P}}(X))_{0}$.
\end{proposition}
\begin{proof}
Assume that $T\in (Q_{\mathcal{P}}(X))_{0}$. It follows from proposition (\ref{proposition:jo2}) that there is
some calibration $\mathcal{P}^{\prime }\in \mathcal{C}(X)$ such that $%
\mathcal{P}\approx \mathcal{P}^{\prime }$ and $T\in B_{\mathcal{P}^{\prime
}}(X)$. Moreover, $Q_{\mathcal{P}}(X)=Q_{\mathcal{P}^{\prime }}(X)$.

If $\mid \lambda \mid >2\left\Vert T\right\Vert _{\mathcal{P}^{\prime }}$,
then the Neumann series $\sum\limits_{n=0}^{\infty }\frac{{T}^{n}}{\lambda
^{n+1}}$ converges in $B_{\mathcal{P}^{\prime }}(X)$ and its sum is $R\left(
\lambda ,T\right) $. This means that the operator $\lambda I-T$ is
invertible in $Q_{\mathcal{P}}(X)$ for all $\mid \lambda \mid >2\left\Vert
T~\right\Vert _{\mathcal{P}^{\prime }}$. Moreover, for each $\epsilon >0$
there exists an index $\ n_{\epsilon }\in \mathbb{N}$ such that 
\begin{equation*}
\left\Vert R\left( \lambda ,T\right) -\sum\limits_{k=0}^{n}\frac{{T}^{k}}{%
\lambda ^{k+1}}\right\Vert _{\mathcal{P}^{\prime }}\mathtt{<}\epsilon
,\left( \forall \right) n\geq n_{\epsilon },
\end{equation*}%
which implies that for each $n\geq n_{\epsilon }$ we have 
\begin{equation*}
\left\Vert R\left( \lambda ,T\right) \right\Vert _{\mathcal{P}^{\prime
}}\leq \left\Vert R\left( \lambda ,T\right) -\sum\limits_{k=0}^{n_{\epsilon
}}\frac{{T}^{k}}{\lambda ^{k+1}}\right\Vert _{\mathcal{P}^{\prime
}}+\left\Vert \sum\limits_{k=0}^{n_{\epsilon }}\frac{{\ T}^{k}}{\lambda
^{k+1}}\right\Vert _{P^{\prime }}\mathtt{<}
\end{equation*}%
\begin{equation*}
\mathtt{<}\epsilon +\mid \lambda \mid ^{-1}\sum\limits_{k=0}^{n_{\epsilon
}}\left\Vert \frac{{\ T}^{k}}{\lambda ^{k}}\right\Vert _{\mathcal{P}%
^{\prime }}\mathtt{<}\epsilon +(2\left\Vert T\right\Vert _{\mathcal{P}%
^{\prime }})^{-1}\sum\limits_{k=0}^{n_{\epsilon }}2^{-k}<\epsilon
+(\left\Vert T\right\Vert _{\mathcal{P}^{\prime }})^{-1}.
\end{equation*}%
Since $\epsilon >0$ is arbitrarily chosen, we have that 
\begin{equation*}
\left\Vert R\left( \lambda ,T\right) \right\Vert _{\mathcal{P}^{\prime }}%
\mathtt{<}(\left\Vert T\right\Vert _{\mathcal{P}^{\prime }})^{-1},\left(
\forall \right) \mid \lambda \mid >2\left\Vert T\right\Vert _{\mathcal{P}%
^{\prime }}
\end{equation*}

From definition of norm $\left\Vert \;\right\Vert _{\mathcal{P}'}$ it follows
that 
\begin{equation*}
\hat{p^{\prime}}\left( R\left( \lambda ,T\right) \right) \mathtt{<}(\left\Vert
T\right\Vert _{\mathcal{P}^{\prime }})^{-1},
\end{equation*}%
for any $p\in \mathcal{P}^{\prime }$ and for each $\mid \lambda \mid
>2\left\Vert T\right\Vert _{\mathcal{P}^{\prime }}$, which means that the
set 
\begin{equation*}
\{R\left( \lambda ,T\right) |\;\mid \lambda \mid >2\left\Vert T\right\Vert _{%
\mathcal{P}^{\prime }}\}
\end{equation*}%
is bounded in $Q_{\mathcal{P}}(X)=Q_{\mathcal{P}'}(X)$. Therefore, $T$ is regular.

Now suppose that $T\in Q_{\mathcal{P}}(X)$ \ is regular, but it is not
bounded in $Q_{\mathcal{P}}(X)$. By propositions \ref{proposition:st1} and \ref{proposition:raza spectrala} this means that 
\begin{equation*}
\mid \sigma _{W}(Q_{\mathcal{P}},T)\mid =\mid \sigma (Q_{\mathcal{P}},T)\mid =r_{\mathcal{P}}(T)=\infty,
\end{equation*}
which contradicts the supposition we have made. Therefore, T is a bounded
element of $Q_{\mathcal{P}}(X)$.
\end{proof}

\begin{proposition}{\label{proposition:qb2}}
Let $(X,\mathcal{P})$ be a locally convex space. If 
$T\in (Q_{\mathcal{P}}(X))_{0}$, then 
\begin{equation*}
\rho _{W}(Q_{\mathcal{P}},T)=\rho (Q_{\mathcal{P}}^{0},T).
\end{equation*}
\end{proposition}
\begin{proof}
If $\lambda _{0}\in \rho (Q_{\mathcal{P}}^{0},T)$ then from previous
proposition results that $R(\lambda _{0},T)$ is a regular element of the
algebra $Q_{\mathcal{P}}(X)$, so there exists $t>0$ for which the condition
(1) and (2) of the definition \ref{definition:reg} are fullfilled. Those
conditions are equivalent with

\begin{enumerate}
\item[1')] $(\lambda -\lambda _{0})^{-1}I-R(\lambda _{0},T)$ is invertible in $%
Q_{\mathcal{P}}(X)$ for all $|\lambda -\lambda _{0}|<t^{-1},\lambda \neq
\lambda _{0}$;

\item[ 2')] the set 
\begin{equation*}
\{R((\lambda -\lambda _{0})^{-1},R(\lambda _{0},T))|~|\lambda -\lambda
_{0}|<t^{-1},\lambda \neq \lambda _{0}\}
\end{equation*}
\end{enumerate}
is bounded in $Q_{\mathcal{P}}(X)$.

From the condition (2') and lemma \ref{lemma:qb1} results that
the set 
\begin{equation*}
\{(\lambda -\lambda _{0})^{-1}R((\lambda -\lambda _{0})^{-1},R(\lambda
_{0},T))|~|\lambda -\lambda _{0}|<t^{-1},\lambda \neq \lambda _{0}\}
\end{equation*}%
is bounded in $Q_{\mathcal{P}}(X)$. Moreover, each seminorm $\hat{p}%
(T),(p\in \mathcal{P})$, is submultiplicative, so the set 
\begin{equation*}
\{(\lambda -\lambda _{0})^{-1}R(\lambda _{0},T)R((\lambda -\lambda
_{0})^{-1},R(\lambda _{0},T))|~|\lambda -\lambda _{0}|<t^{-1},\lambda \neq
\lambda _{0}\}
\end{equation*}%
is also bounded in $Q_{\mathcal{P}}(X)$. Since 
\begin{equation*}
(\lambda I-T)(\lambda _{0}-\lambda )^{-1}R(\lambda _{0},T)R((\lambda
-\lambda _{0})^{-1},R(\lambda _{0},T))=
\end{equation*}%
\begin{equation*}
=((\lambda _{0}I-T)+(\lambda -\lambda _{0})I)(\lambda _{0}-\lambda
)^{-1}R(\lambda _{0},T)R((\lambda -\lambda _{0})^{-1},R(\lambda _{0},T))=
\end{equation*}%
\begin{equation*}
=(\lambda _{0}-\lambda )^{-1}R((\lambda -\lambda _{0})^{-1},R(\lambda
_{0},T))-R(\lambda _{0},T)R((\lambda -\lambda _{0})^{-1},R(\lambda _{0},T))=
\end{equation*}%
\begin{equation*}
=((\lambda _{0}-\lambda )^{-1}I-R(\lambda _{0},T))R((\lambda -\lambda
_{0})^{-1},R(\lambda _{0},T))=I,
\end{equation*}%
results that 
\begin{equation}
{\label{equation:qb2}}R(\lambda ,T)=(\lambda _{0}-\lambda )^{-1}R(\lambda
_{0},T)R((\lambda -\lambda _{0})^{-1},R(\lambda _{0},T)).
\end{equation}%
Therefore, the conditions \newline
1) $\lambda I-T$ is invertible for all $|\lambda -\lambda _{0}|<t^{-1}$;%
\newline
2) $\{R(\lambda ,T)|~|\lambda -\lambda _{0}|<t^{-1}\}$ is bounded in $Q_{%
\mathcal{P}}(X)$,\newline
of definition (\ref{definition:sw}) are fullfilled, so $\lambda _{0}\in \rho
_{W}(Q_{\mathcal{P}},T)$ and $\rho (Q_{\mathcal{P}}^{0},T)\subset \rho
_{W}(Q_{\mathcal{P}},T)$.

Conversely, if $\lambda _{0}\in \rho _{W}(Q_{\mathcal{P}},T)$ there exists $%
K>0$ such that\newline
1'') $\lambda I-T$ is invertible for all $|\lambda -\lambda _{0}|<K$;\newline
2'') $\{R(\lambda ,T)|~|\lambda -\lambda _{0}|<K\}$ is bounded in $Q_{%
\mathcal{P}}(X)$.

From the relation (\ref{equation:qb2}) and property (2'') will result that
the 
\begin{equation*}
\{(\lambda -\lambda _{0})^{-1}R((\lambda -\lambda _{0})^{-1},R(\lambda
_{0},T))|~|\lambda -\lambda _{0}|^{-1}>K^{-1},\lambda \neq \lambda _{0}\}
\end{equation*}%
is bounded in $Q_{\mathcal{P}}(X)$, so $R(\lambda _{0},T)$ is regular in $Q_{%
\mathcal{P}}(X)$. From previous proposition results that $R(\lambda
_{0},T)\in (Q_{\mathcal{P}}(X))_{0}$ and $\lambda _{0}\in \rho (Q_{\mathcal{P%
}}^{0},T)$.
\end{proof}

\begin{proposition}{\label{proposition:qb3}}
Let $(X,\mathcal{P})$ be a locally convex space. If $T\in (Q_{\mathcal{P}}(X))_{0}$ and $\mid \lambda
_{0}\mid >r_{\mathcal{P}}(T)$, then $\lambda _{0}\in \rho (Q_{\mathcal{P}%
}^{0},T)$.
\end{proposition}
\begin{proof}
Proposition \ref{proposition:newman} implies that the series $%
\sum\limits_{n=0}^{\infty }\frac{{\ T}^{n}}{\lambda _{0}^{n+1}}$ converges
to $R\left( \lambda _{0},T\right) \in Q_{\mathcal{P}}(X)$. Then, there
exists $\epsilon \mathtt{>}0$ such that 
\begin{equation*}
D\left( \lambda _{0},\varepsilon \right) =\left\{ \lambda ||\lambda -\lambda
_{0}|<\varepsilon \right\} \subset \left\{ \mu ||\mu |>r_{\mathcal{P}%
}(T)\right\} ,
\end{equation*}%
so the operator $\lambda I-T$ is invertible, for every $\lambda \in D\left( \lambda _{0},\epsilon \right) $, and $(\lambda I-T)^{-1}\in Q_{\mathcal{P%
}}(X)$.

Now we will prove that the set $\sigma (Q_{\mathcal{P}},R(\lambda _{0},T))$
is bounded. If $\mid \mu \mid >\epsilon ^{-1}$, then $\mid \mu \mid
^{-1}<\epsilon $ and $\lambda _{0}-\mu ^{-1}\in D\left( \lambda
_{0},\epsilon \right) $. From previous observations results that $(\lambda
_{0}-\mu ^{-1})I-T$ is invertible and $((\lambda _{0}-\mu ^{-1})I-T)^{-1}\in
Q_{\mathcal{P}}(X)$.

Since 
\begin{equation*}
\mu ^{-1}R\left( \lambda _{0}-\mu ^{-1},T\right) \left( \lambda
_{0}I-T\right) (\mu I-R\left( \lambda _{0},T\right) )=
\end{equation*}%
\begin{equation*}
=R\left( \lambda _{0}-\mu ^{-1},T\right) \left( \lambda _{0}I-T\right) -\mu
^{-1}R\left( \lambda _{0}-\mu ^{-1},T\right) =
\end{equation*}%
\begin{equation*}
=R\left( \lambda _{0}-\mu ^{-1},T\right) (((\lambda _{0}-\mu ^{-1})I-T)+\mu
^{-1}I)-\mu ^{-1}R\left( \lambda _{0}-\mu ^{-1},T\right) =
\end{equation*}%
\begin{equation*}
=I+\mu ^{-1}R\left( \lambda _{0}-\mu ^{-1},T\right) -\mu ^{-1}R\left(
\lambda _{0}-\mu ^{-1},T\right) =I.
\end{equation*}%
results that%
\begin{equation*}
R(\mu ,R\left( \lambda _{0},T\right) )=\mu ^{-1}R\left( \lambda _{0}-\mu
^{-1},T\right) \left( \lambda _{0}I-T\right)
\end{equation*}

But 
\begin{equation*}
R\left( \lambda _{0}-\mu ^{-1},T\right) ,\left( \lambda _{0}I-T\right)
\in Q_{\mathcal{P}}(X),
\end{equation*}
so $R(\mu ,R\left( \lambda _{0},T\right) )\in Q_{%
\mathcal{P}}(X)$, for all $\mid \mu \mid >\epsilon ^{-1}$. 

Therefore, $\sigma (Q_{\mathcal{P}},R(\lambda _{0},T))\subset D\left( 0,\epsilon ^{-1}\right) $
and from proposition \ref{proposition:gi2} results that there $\mathcal{P}^{\prime }\in \mathcal{C}(X)$ such that $R(\lambda _{0},T)\in B_{\mathcal{P}^{\prime
}}(X)$ and $\mathcal{P}^{\prime }\approx \mathcal{P}$. This implies that $Q_{%
\mathcal{P}^{\prime }}(X)=Q_{\mathcal{P}}(X)$ and $R(\lambda _{0},T)\in (Q_{%
\mathcal{P}}(X))_{0}$ (proposition \ref{proposition:jo2}).
\end{proof}

\begin{corollary}
{\label{corollary:qb 2}}Let $X$ be a locally convex space and $\mathcal{P}%
\in \mathcal{C}(X)$. If $T\in (Q_{\mathcal{P}}(X))_{0}$ then 
\begin{equation*}
\left\vert \sigma (Q_{\mathcal{P}},T)\right\vert =\left\vert \sigma _{W}(Q_{%
\mathcal{P}},T)\right\vert =r_{\mathcal{P}}(T)
\end{equation*}
\end{corollary}
\begin{proof}
Is a direct consequence of \ propositions \ref{proposition:qb2}, \ref%
{proposition:qb3} and \ref{proposition:raza spectrala}
\end{proof}

\begin{definition}
Let $(X,\mathcal{P})$ be a locally convex space. An operator  $T\in Q_{%
\mathcal{P}}(X)$ is $\mathcal{P}$-quasnilpotent if  $r_{\mathcal{P}}(T)=0$.
\end{definition}

\begin{remark}
\begin{enumerate}
\item 
If $T\in Q_{\mathcal{P}}(X)$ is $\mathcal{P}$-quasnilpotent, then $%
T\in (Q_{\mathcal{P}}(X))_{0}$ and $\sigma _{W}(Q_{\mathcal{P}},T)=\{0\}$.
\item 
$T\in Q_{\mathcal{P}}(X)$ is $\mathcal{P}$-quasnilpotent if and
only if  $\sigma (Q_{\mathcal{P}},T)=\{0\}$.
\end{enumerate}
\end{remark}

\section{A functional calculus}

L.Waelbroeck presented before a functional calculus for regular
operator on quasi-complet locally convex space. In this section we prove
that it can be developed a functional calculus for the bounded elements of
algebra $Q_{\mathcal{P}}(X)$ (which by previous section are regular element
of this algebra), when $(X,\mathcal{P})$ ($\mathcal{P}\in \mathcal{C}(X)$) is a sequentially complete locally convex algebra. We will use some techniques from I. Colojoara \cite{co} and
L.Waelbroeck \cite{wa}.

Let $\ \mathcal{P}\in C(X)$ be arbitrary chosen and $D\subset \mathbb{C}$ a
relatively compact open set. Denote by $\mathcal{O}(D,Q_{\mathcal{P}}(X))$
the unitary algebra of the functions $f:D\rightarrow Q_{\mathcal{P}}(X)$
which are holomorphic on $D$ and continuous on $\overline{D}$.

\begin{lemma}
If $p\in \mathcal{P}$, then the application $\mid f\mid _{p,D}:\mathcal{O}%
(D,Q_{\mathcal{P}}(X))\rightarrow R$ given by relation, 
\begin{equation*}
\mid f\mid _{p,D}=\sup\limits_{z\in D}p(f\left( z\right) ),\left( \forall
\right) f\in \mathcal{O}(D,Q_{\mathcal{P}}(X)),
\end{equation*}
is a submultiplicative seminorm on $\mathcal{O}(D,Q_{\mathcal{P}}(X)).$
\end{lemma}

If we denote by $\tau _{\mathcal{P},D}$ the topology defined by the family $%
\{$ $\mid f\mid _{p,D}|$ $p\in \mathcal{P}\}$ on $\mathcal{O}(D,Q_{\mathcal{P%
}}(X))$, then $(\mathcal{O}(D,Q_{\mathcal{P}}(X)),\tau _{\mathcal{P},D})$ is
a l.m.c.-algebra.

Let $K\subset \mathbb{C}$ be a compact set, arbitrary chosen. We define the set
\begin{center}
$\mathcal{O}(K,Q_{\mathcal{P}}(X))=\cup \{\mathcal{O}(D,Q_{\mathcal{P}}(X))|$
$D\subset\mathbf{C}$ is relatively compact open set $\}$
\end{center}

If $D_{1},D_{2}\subset \mathbf{C}$ are relatively compact open sets such
that $K\subset D_{i},$ $i=\overline{1,2}$, and $f_{i}\in \mathcal{O}%
(D_{i},Q_{\mathcal{P}}(X))$, $i=\overline{1,2}$, we say that $f_{1}\backsim
f_{2}$ if and only if there exists an open set D such that $K\subset
D\subset D_{1}\cap D_{2}$ and $f_{1}|_{D}=f_{2}|_{D}$. Let denote by $%
\mathcal{A}(K,Q_{\mathcal{P}}(X))$ be the set of the equivalence classes of $%
\mathcal{O}(K,Q_{\mathcal{P}}(X))$ in respect with this equivalence
relation. It is easily to see that $\mathcal{A}(K,Q_{\mathcal{P}}(X))$ is a
unitary algebra and the elements of this algebra are usually called germs of
the holomorphic functions from $K$ to $Q_{\mathcal{P}}(X)$.

\begin{remark}
We consider the following notations:
\begin{enumerate}
\item 
$\tilde{f}$ is the germ of the holomorphic function $f\in \mathcal{O}%
(D,Q_{\mathcal{P}}(X))$.
\item 
$\varphi $ is the canonical morphism $\mathcal{O}(K,Q_{\mathcal{P}%
}(X))\rightarrow \mathcal{A}(K,Q_{\mathcal{P}}(X))$;
\item 
$\varphi _{D}$ is the restriction of $\varphi $\ to $\mathcal{O}(D,Q_{%
\mathcal{P}}(X))$.
\end{enumerate}
\end{remark}

\begin{remark}
\begin{enumerate}
\item Since we can identifies $\mathbb{C}$ with $\mathbb{C}I=\{\lambda
I~|~\lambda \in \mathbb{C}~\}$, the algebras $\mathcal{O}(K,\mathbb{C})$ and $%
\mathcal{A}(K,\mathbb{C})$ can be considerate subalgebras of $\mathcal{O}%
(K,Q_{\mathcal{P}}(X))$, respectively $\mathcal{A}(K,Q_{\mathcal{P}}(X))$.
Therefore, we write $\mathcal{O}(K)$ and $\mathcal{A}(K)$ instead of $%
\mathcal{O}(K,\mathbb{C})$ and $\mathcal{A}(K,\mathbb{C})$

\item If $\tau _{\mathcal{P},ind}=\underset{\rightarrow D}{\lim }\tau _{%
\mathcal{P},D}$ (inductive limit), then $(\mathcal{A}(K,Q_{\mathcal{P}%
}(X)),\tau _{\mathcal{P~}ind})$ is a l.m.c.-algebra.
\end{enumerate}
\end{remark}

We need the following lemma from complex analysis.

\begin{lemma}
For each compact set $K\subset\mathbb{C}$ and each relatively compact open
set $D\supset K$ there exists some open set $G$ such that:

\begin{enumerate}
\item $K\subset G\subset \overline{G}\subset D$;

\item G has a finite number of conex components $(G_{i})_{i=\overline{1,n}}$%
, the closure of which are pairwise disjoint;

\item the boundary $\partial G_{i}$ of $G_{i}, i=\overline{1,n}$, consists
of a finite positive number of closed rectifiable Jordan curves $(\Gamma
_{ij})_{j=\overline{1,m_{i}}}$, no two of which intersect;

\item $K\cap \Gamma _{ij}=\varnothing  $, for each $i=\overline{1,n}$ and every $j=%
\overline{1,m_{i}}$.
\end{enumerate}
\end{lemma}

\begin{definition}
If the sets $K$ and $D$ are like in the previous lemma, then an open set $G$
is called Cauchy domain for pair $(K,D)$ if it satisfies the properties
(1)-(4). The boundary 
\begin{equation*}
\Gamma =\cup _{i=\overline{1,n}}\cup _{j=\overline{1,m_{i}}}\Gamma _{ij}
\end{equation*}
of $G$ is called Cauchy boundary for pair $(K,D)$.
\end{definition}

\begin{theorem}
If $\ \mathcal{P}\in \mathcal{C}_{0}(X)$ and $T\in (Q_{\mathcal{P}}(X))_{0}$%
, then for each relatively compact open set $D\supset \sigma _{W}(Q_{%
\mathcal{P}},T)$ there exists an application 
\begin{equation*}
\mathcal{F}_{T,D}:\mathcal{O}(D,Q_{\mathcal{P}}(X))\rightarrow Q_{\mathcal{P}%
}(X)
\end{equation*}
with the properties:
\begin{enumerate}
\item 
$\mathcal{F}_{T,D}$ is continuous and linear;
\item 
$\mathcal{F}_{T,D}\left( k_{S}\right) =S$, where $k_{S}\equiv S$;
\item 
$\mathcal{F}_{T,D}\left( id_{I}\right) =T$, where $id_{I}(\lambda
)=\lambda I$, for every $\lambda \in \mathbb{C}$.
\end{enumerate}
\end{theorem}
\begin{proof}
Let $\Gamma $ be a Cauchy boundary for the pair $(\sigma _{W}(Q_{\mathcal{P}%
},T),D)$. Then the integral%
\begin{equation*}
\frac{1}{2\pi i}\int\limits_{\Gamma }f\left( \lambda \right) R\left( \lambda
,T\right) d\lambda ,\left( \forall \right) f\in \mathcal{O}(D,Q_{\mathcal{P}}(X)),
\end{equation*}%
exists like Stieltjes integral, since $Q_{\mathcal{P}}(X)$ is a sequentially
complete l.m.c.-algebra and the applications $t\rightsquigarrow f(\omega
(t))R(\omega (t),T)$ are continuous on $[0,1]$ for a continuous
parametrization $\omega $ of $\Gamma $.

Moreover, if $\ \Gamma _{1}$ and $\Gamma _{2}$ are Cauchy boundaries for
pair $(\sigma _{W}(Q_{\mathcal{P}},T),D)$ then 
\begin{equation*}
\frac{1}{2\pi i}\int\limits_{\Gamma _{1}}f\left( z\right) R\left( \lambda
,T\right) d\lambda =\frac{1}{2\pi i}\int\limits_{\Gamma _{2}}f\left( \lambda
\right) R\left( \lambda ,T\right) d\lambda ,\left( \forall \right) f\in 
\mathcal{O}(D,Q_{\mathcal{P}}(X)),
\end{equation*}
therefore the application $\mathcal{F}_{T,D}:\mathcal{O}(K,Q_{\mathcal{P}%
}(X))\rightarrow Q_{\mathcal{P}}(X)$ given by formula 
\begin{equation*}
\mathcal{F}_{T,D}(f)=\frac{1}{2\pi i}\int\limits_{\Gamma }f\left( \lambda
\right) R\left( \lambda ,T\right) dz,\left( \forall \right) f\in \mathcal{O}%
(D,Q_{\mathcal{P}}(X)),
\end{equation*}
is well defined. Now we prove that $\mathcal{F}_{T,D}$ has the properties
(1)-(3).

The linearity is obvious. For every $p\in \mathcal{P}$\ and \ every $f\in 
\mathcal{O}(D,Q_{\mathcal{P}}(X))$ we have 
\begin{equation*}
\hat{p}(\mathcal{F}_{T,D}(f))\leq \frac{L(\Gamma )}{2\pi }\underset{\lambda
\in \Gamma }{\sup }\hat{p}(R\left( \lambda ,T\right) )\underset{\lambda \in
\Gamma }{\sup }\hat{p}(f(\lambda ))\leq \frac{L(\Gamma )}{2\pi }\underset{%
\lambda \in \Gamma }{\sup }\hat{p}(R\left( \lambda ,T\right) )\mid f\mid
_{p,D},
\end{equation*}%
where $L(\Gamma )$ is the lenght of $\Gamma $, which implies the continuity
of application $\mathcal{F}_{T,D}$.

Let $r>r_{\mathcal{P}}(T)$ and $\Gamma _{r}=\{z\in C|~|z|=r\}$. For each $%
\lambda \in \Gamma _{r}$ we have $r_{\mathcal{P}}(\frac{T}{\lambda })<1$, so
from lemma \ref{lemma:qb0} results that 
\begin{equation*}
R\left( \lambda ,T\right) =\lambda ^{-1}(I-\frac{T}{\lambda })=\lambda ^{-1}%
\underset{n\in \mathbb{N}}{\sum }\left( \frac{T}{\lambda }\right) ^{n}=%
\underset{n\in \mathbb{N}}{\sum }\frac{T^{n}}{\lambda ^{n+1}}
\end{equation*}

This observation implies that
\begin{equation*}
\mathcal{F}_{T,D}\left( k_{S}\right) =\frac{1}{2\pi i}\int\limits_{\Gamma
}k_{S}\left( \lambda \right) R\left( \lambda ,T\right) d\lambda =\frac{S}{%
2\pi i}\underset{n\in \mathbb{N}}\sum T^{n}\int\limits_{\Gamma }\frac{d\lambda}{%
\lambda ^{n+1}}=S
\end{equation*}
\begin{equation*}
\mathcal{F}_{T,D}\left( id_{I}\right) =\frac{1}{2\pi i}\int\limits_{\Gamma
}\lambda R\left( \lambda ,T\right) d\lambda =\frac{1}{2\pi i}\underset{n\in 
\mathbb{N} }\sum T^{n}\int\limits_{\Gamma }\frac{d\lambda}{\lambda ^{n}}=T.
\end{equation*}
\end{proof}

\begin{corollary}
If $\mathcal{P}\in \mathcal{C}_{0}(X)$ and $T\in (Q_{\mathcal{P}}(X))_{0}$,
then there exists an application $\mathcal{F}_{T}:\mathcal{A}(\sigma _{W}(Q_{%
\mathcal{P}},T),Q_{\mathcal{P}}(X))\rightarrow Q_{\mathcal{P}}(X)$ which satisfies the
conditions:
\begin{enumerate}
\item 
$\mathcal{F}_{T}$ is continuous and linear;
\item 
$\mathcal{F}_{T}\left( \tilde{k}_{S}\right) =S$, where $\tilde{k}_{T}$
is the germ of the function $k_{S}\equiv S$;
\item 
$\mathcal{F}_{T}\left( \widetilde{id}_{I}\right) =T$, where $\widetilde{id}%
_{I} $ is the germ of the function $id_{I}(\lambda )=\lambda I$, for all $\lambda\in \mathbb{C}$
\end{enumerate}
\end{corollary}
\begin{proof}
If $\tilde{f}\in \mathcal{A}(\sigma _{W}(Q_{\mathcal{P}},T))$, then we
consider 
\begin{equation*}
\mathcal{F}_{T}(\tilde{f})=\mathcal{F}_{T,D}\left( f\right) ,\left( \forall
\right) \tilde{f}\in \mathcal{A}(\sigma _{W}(Q_{\mathcal{P}},T)),
\end{equation*}%
where $f\in \mathcal{O}(D,Q_{\mathcal{P}},T))$ is an element of equivalence class $\tilde{f}$.
It is obvious that the definition of $\mathcal{F}_{T,D}\left( f\right) $ is
independent by the function $f$ and $\mathcal{F}_{T}\left( \tilde{f}\right) $ is
linear. Since $\mathcal{F}_{T,D}=\mathcal{F}_{T}\circ \varphi _{D}$ and $%
\mathcal{F}_{T,D}$ is continuous results that $\mathcal{F}_{T}$ is
continuous.

The properties (2) and (3) results directly from the previous theorem.
\end{proof}

\begin{corollary}
{\label{corollary:fc}}If $\mathcal{P}\in \mathcal{C}_{0}(X)$ and $T\in (Q_{%
\mathcal{P}}(X))_{0}$, then there exists an unique unitary continuous
morphism $F_{T}:\mathcal{A}(\sigma _{W}(Q_{\mathcal{P}},T))\rightarrow Q_{%
\mathcal{P}}(X)$ which satisfies the condition $F_{T}\left( \widetilde{
id}\right) =T$, where $id$ is the identity function on $\mathbb{C}$.
\end{corollary}
\begin{proof}
The application $F_{T}$ and $F_{T,D}$ are defined in the same way like the
application $\mathcal{F}_{T}$\ and $\mathcal{F}_{T,D}$. It is easily to see
that $F_{T}$\ and $F_{T,D}$ are linear and continuous.
Moreover, $F_{T}$ is unitary and $F_{T}\left( \widetilde{id}\right) =T.$

We prove that $\mathcal{F}_{T}$ is multiplicative. Let $\tilde{f},\tilde{g}%
\in \mathcal{A}(\sigma _{W}(Q_{\mathcal{P}},T))$ and $f\in \tilde{f} $
respectively $g\in \tilde{g}$. We consider that $G$ and $G^{\prime }$ are
two Cauchy domains with the property $\overline{G^{\prime }}\subset G$. If $%
\Gamma $ and $\Gamma ^{\prime }$ are the boundaries of $G$ and $G^{\prime }$
then 
\begin{equation*}
F_{T}( \tilde{f})F_{T}\left( \tilde{g}\right) =-\frac{1}{(2\pi i)^{2}}\int_{\Gamma
}\int_{\Gamma ^{\prime }}f(\lambda )g(\omega )R(\lambda ,T)R(\omega
,T)d\lambda d\omega
\end{equation*}

Since $\overline{G^{\prime }}\subset G$, results that $\Gamma\cap \Gamma
^{\prime }=\Phi $, so 
\begin{equation*}
\omega -\lambda \neq 0,~(\forall )\lambda \in \Gamma, (\forall )\omega \in
\Gamma ^{\prime }.
\end{equation*}

Therefore using the equality 
\begin{equation*}
R(\lambda ,T)-R(\omega ,T)=(\omega -\lambda )R(\lambda ,T)R(\omega ,T)
\end{equation*}
we have 
\begin{equation*}
F_{T}\left( \tilde{f}\right) F_{T}\left( \tilde{g}\right) =\frac{1}{(2\pi
i)^{2}}\int_{\Gamma }f(\lambda )R(\lambda ,T)\left( \int_{\Gamma ^{\prime }}%
\frac{g(\omega )}{\omega -\lambda }d\omega \right) d\lambda+
\end{equation*}
\begin{equation*}
+\frac{1}{(2\pi i)^{2}}\int_{\Gamma ^{\prime }}g(\omega )R(\omega ,T)\left(
\int_{\Gamma }\frac{f(\lambda )}{\lambda -\omega }d\lambda \right) d\omega =
\end{equation*}
\begin{equation*}
=\frac{1}{2\pi i}\int_{\Gamma }f(\lambda )g\left( \lambda
\right) R(\lambda ,T)d\omega =F_{T}\left( \tilde{f}\tilde{g}\right)
\end{equation*}

Assume that $F:\mathcal{A}(\sigma _{W}(Q_{\mathcal{P}},T))\rightarrow Q_{%
\mathcal{P}}(X)$ is an unitary continuous morphism which satisfies the
condition $F\left( \widetilde{id}\right) =T$. We prove that $F_{T}=F$.

Let $\tilde{f}\in \mathcal{A}(\sigma _{W}(Q_{\mathcal{P}},T))$, $D\supset
\sigma _{W}(Q_{\mathcal{P}},T)$ a relatively compact open set, $f\in 
\mathcal{O}(D)$, such that $f\in\tilde{f}$, and G a Cauchy domain for $%
(\sigma _{W}(Q_{\mathcal{P}},T),D)$ with the boundary $\Gamma $. For every $%
n\in N^{\ast }$ and $z_{1},...,z_{n}\in \Gamma $ we consider the function $%
f_{n}:G\rightarrow \mathbb{C}$ given by the relation 
\begin{equation}
{\label{equation:fc3}} f_{n}\left( \omega \right) =\frac{1}{2\pi i}%
\sum_{j=1}^{n}\frac{f\left( z_{j}\right) \left( z_{j+1}-z_{j}\right) }{%
z_{j}-\omega } , \left( \forall \right) \omega \in G.
\end{equation}

Then, 
\begin{equation*}
\underset{n\rightarrow \infty }\lim f_{n}\left( \omega \right) =\frac{1}{%
2\pi i}\int_{\Gamma }\frac{f(z)}{z-\omega }dz=f\left( \omega \right)
\end{equation*}
and since the convergence is uniformly on each compact set $K\subset G$,
results that $\lim_{\tau _{ind}}\tilde{f}_{n}=\tilde{f}$. Using the continuity of $F$ results that 
\begin{equation}
{\label{equation:fc1}} \underset{n\rightarrow \infty }{\lim }F\left( \tilde{f%
}_{n}\right) =F\left( \tilde{f}\right)
\end{equation}

Since $F$ is a unitary morphism with the property $F\left( \widetilde{id}%
\right) =T$, then from relation (\ref{equation:fc3}) results that 
\begin{equation*}
F\left( \tilde{f}_{n}\right) =\frac{1}{2\pi i}\ \sum_{j=1}^{n}f\left(
z_{j}\right) \left( z_{j+1}-z_{j}\right) R(z_{j},T)
\end{equation*}
so 
\begin{equation}
{\label{equation:fc2}} \underset{n\rightarrow \infty }\lim F\left( \tilde{f}%
_{n}\right) =\frac{1}{2\pi i}\int\limits_{\Gamma }f\left( z\right) R\left(
z,T\right) dz,
\end{equation}

From relations (\ref{equation:fc1}) and (\ref{equation:fc2}) results that 
\begin{equation*}
F\left( \tilde{f}\right) =\frac{1}{2\pi i}\int\limits_{\Gamma }f\left(
z\right) R\left( z,T\right) dz=F_{T,D}(f)=F_{T}\left( \tilde{f}\right)
\end{equation*}
which implies that $F_{T}=F$.
\end{proof}

\begin{lemma}
{\label{lemma:reg1}} If $K\subset \mathbb{C}$ is a compact set, then each
element of the algebra $\mathcal{A}(K)$ is regular.
\end{lemma}
\begin{proof}
Let $\tilde{f}\in \mathcal{A}(K)$, $D\supset K$ a relatively compact open
set, $f\in \mathcal{O}(D)$ ($f\in $ $\tilde{f}$) and $\omega _{0}\notin
f\left( K\right)$. Then there exists two relatively compact open set $U$
and $V$ such that $\omega _{0}\in U$, $f\left( K\right) \subset V$ and $%
\overline{U}\cap \overline{V}=\Phi $. For every $\omega \in U$ the function $%
f_{\omega }:\overline{f^{-1}(V)}\rightarrow \mathbb{C}$ given by relation 
\begin{equation*}
f_{\omega }\left( \lambda \right) =\frac{1}{\omega -f\left( \lambda \right) }%
,\left( \forall \right) \lambda \in \overline{f^{-1}(V)}
\end{equation*}%
is holomorphic on $\ f^{-1}(V)$, so $\tilde{f}_{\omega }\in \mathcal{A}(K)$.

Since for every compact set $A\subset f^{-1}(V)$ we have 
\begin{equation*}
\underset{\omega \in \overset{\_\_}{U}}{\sup }\underset{\lambda \in A}{~\sup 
}|f_{\omega }\left( \lambda \right) |<\infty
\end{equation*}%
results that the set $\{$ $\tilde{f}_{\omega }|\omega \in \overline{U}\}$ is
bounded in $(\mathcal{A}(K),\tau _{ind})$. Moreover, 
\begin{equation*}
(\omega \tilde{1}-\tilde{f})\tilde{f}_{\omega }=\overset{\sim }{1}
\end{equation*}%
so $\omega \in \sigma _{W}(\tilde{f})$. Therefore $\sigma _{W}(\tilde{f}%
)\subset f\left( K\right) $. Since $K$ is compact the set $f\left( K\right) $
is compact, so $\sigma _{W}(\tilde{f})$ is compact and $\tilde{f}$ is
regular.
\end{proof}

\begin{lemma}
{\label{lemma:reg2}} If $X$ and $Y$ are unitary locally convex algebra and $%
F:X\rightarrow Y$ is unitary continuous morphism, then $F\left( X_{r}\right)
\subset Y_{r}$, where $X_{r}$ and $Y_{r}$ are the algebras of the regular elements of $X$, respectively $Y$.
\end{lemma}
\begin{proof}
If $x\in X_{r}$, then there exists $k>0$ such that $\lambda e-x$ is
invertible for every $|\lambda |>k$ and the set $\{R\left( \lambda ,x\right)
||\lambda |>k\}$ is bounded in $X$. Since $F$ is unitary morphism results
that 
\begin{equation*}
F\left( R\left( \lambda ,x\right) \right) =R\left( \lambda ,F(x\right)
), \left( \forall \right) |\lambda |>k, 
\end{equation*}
so from continuity of $F$ results that the set 
\begin{equation*}
\{F\left( R\left( \lambda ,x\right) \right) ||\lambda |>k\}=\{R\left(
\lambda ,F(x\right) )||\lambda |>k\}
\end{equation*}%
is bounded. Therefore $F\left( x\right) $ is regular.
\end{proof}

\begin{proposition}
{\label{proposition:reg1}}If $\mathcal{P}\in \mathcal{C}_{0}(X)$ and $T\in
(Q_{\mathcal{P}}(X))_{0}$, then 
\begin{equation*}
F_{T}\left( \mathcal{A}(\sigma _{W}(Q_{\mathcal{P}},T))\right) \subset (Q_{%
\mathcal{P}}(X))_{0}.
\end{equation*}
\end{proposition}
\begin{proof}
From lemmas \ref{lemma:reg1} and \ref{lemma:reg2} results that $F_{T}(\tilde{%
f})$ is a regular element of algebra $Q_{\mathcal{P}}(X)$, for every $\tilde{%
f}\in \mathcal{A}(\sigma _{W}(Q_{\mathcal{P}},T))$, so by proposition \ref%
{proposition:qb1} we have that $\tilde{f}\in (Q_{\mathcal{P}}(X))_{0}$.
\end{proof}

\begin{lemma}
If $\mathcal{P}\in \mathcal{C}_{0}(X)$, $T\in (Q_{\mathcal{P}}(X))_{0}$ and $%
P$ is a polynomial, then 
\begin{equation*}
F_{T}(\tilde{P})=\tilde{P}(T) \text{ and } F_{T,D}(P)=P(T).
\end{equation*}
for each relatively compact open set $D\supset \sigma _{W}(Q_{\mathcal{P}},T)$.
\end{lemma}

\begin{remark}
From previous lemma results that for each $T\in (Q_{\mathcal{P}}(X))_{0}$ we
can use the following notation: 
\begin{equation*}
F_{T}(\tilde{f})=\tilde{f}(T)\text{ and }F_{T,D}(f)=f(T).
\end{equation*}%
where $\tilde{f}\in \mathcal{A}(K)$, $D\supset K$ open set and $f\in 
\mathcal{O}(D)$, such that $f\in $ $\tilde{f}$ .
\end{remark}

The following theorem represents the analogous of the spectral mapping theorem
for Banach spaces.

\begin{theorem}
If $\mathcal{P}\in \mathcal{C}_{0}(X)$, $T\in (Q_{\mathcal{P}}(X)~)_{0}$ and 
$f$ is a holomorphic function on an open set $D\supset \sigma _{W}(Q_{%
\mathcal{P}},T)$, then 
\begin{equation*}
\sigma _{W}(Q_{\mathcal{P}},f(T))=f(\sigma _{W}(Q_{\mathcal{P}},T)).
\end{equation*}
\end{theorem}
\begin{proof}
From lemma (\ref{lemma:reg2}) results that the operator $\tilde{f}(T)$ is regular
element of the algebra $Q_{\mathcal{P}}(X)$, so the spectrum $\sigma _{W}(Q_{%
\mathcal{P}},f(T))$ is compact. 

Let $\omega _{0}\notin f(\sigma _{W}(Q_{\mathcal{P}},T))$. Then there exists two relatively compact open set $U$ and 
$V$ such that $\omega _{0}\in U$, $\sigma _{W}(Q_{\mathcal{P}},f(T))\subset
V $ and $\overline{U}\cap \overline{V}=\Phi $. We proved already in the
proof of lemma (\ref{lemma:reg1}) that if the functions $f_{\omega }:%
\overline{f^{-1}(V)}\rightarrow \mathbb{C}$, $\omega \in U$, are given by
relation 
\begin{equation*}
f_{\omega }\left( \lambda \right) =\frac{1}{\omega -f\left( \lambda \right) }%
,\left( \forall \right) \lambda \in \overline{f^{-1}(V)}
\end{equation*}%
then the set $\{\tilde{f}_{\omega }|\omega \in U\}$ is bounded in $(\mathcal{%
A}(\sigma _{W}(Q_{\mathcal{P}},f(T))),\tau _{ind})$. The morphism $F_{T}$ is
unitary, so 
\begin{equation*}
F_{T}(\tilde{f}_{\omega })(\omega I-F_{T}(\tilde{f}))=F_{T}(\tilde{1})=I.
\end{equation*}

Now from the continuity of $F_{T}$ results that the set 
\begin{equation*}
\{F_{T}(\tilde{f}_{\omega })~|~\omega \in \overline{U}~\}=\{R(\omega ,F_{T}(%
\tilde{f}))~|~\omega \in \overline{U}~\}=\{R(\omega ,\tilde{f}(T))~|~\omega
\in \overline{U}~\}
\end{equation*}
is bounded in $Q_{\mathcal{P}}(X)$. Therefore, $\omega _{0}\notin \sigma
_{W}(Q_{\mathcal{P}},f(T))$ and 
\begin{equation*}
\sigma _{W}(Q_{\mathcal{P}},f(T))\subset f(\sigma _{W}(Q_{\mathcal{P}},T)).
\end{equation*}
If $\omega _{0}\in \sigma _{W}(Q_{\mathcal{P}},T)$ and $g_{\omega
_{0}}:D\rightarrow \mathbb{C}$ is defined by 
\begin{equation*}
g_{\omega _{0}}(\lambda )=\left \{%
\begin{array}{c}
\frac{f(\lambda )-f(\omega _{0})}{\lambda -\omega _{0}},\text{\ for }\lambda
\neq \omega _{0}, \\ 
f^{\prime }(\omega _{0}),\text{ \ for }\lambda =\omega _{0},%
\end{array}%
\right.
\end{equation*}
then $g_{\omega _{0}}\in \mathcal{O}(D)$ and 
\begin{equation*}
f(\omega _{0})-f(\lambda )=(\omega _{0}-\lambda )g_{\omega _{0}}(\lambda
),\left( \forall \right) \lambda \in D.
\end{equation*}

Therefore 
\begin{equation*}
f(\omega _{0})I-f(T)=(\omega _{0}I-T)g_{\omega _{0}}(T).
\end{equation*}

Since $\omega _{0}I-T$ is not invertible results that $f(\omega _{0})\in
\sigma _{W}(Q_{\mathcal{P}},f(T))$, so
\begin{center}
$f(\sigma _{W}(Q_{\mathcal{P}},T))\subset \sigma _{W}(Q_{\mathcal{P}},f(T)).$
\end{center}
\end{proof}

\begin{theorem}
Let $\mathcal{P}\in \mathcal{C}_{0}(X)$ and $T\in (Q_{\mathcal{P}%
}(X))_{0}$. If $f$ is holomorphic function on the open set $D\supset \sigma
_{W}(Q_{\mathcal{P}},T)$ and $g\in O(D_{g})$, such that $D_{g}\supset f(D)$,
then $\left( g\circ f\right) (T)=g(f(T))$.
\end{theorem}
\begin{proof}
Let G be a Cauchy domain for the pair ($\sigma _{W}(Q_{\mathcal{P}},T),D)$
and $\Gamma $ the boundary of $G$. Since for each $\omega \notin \overline{%
f(G)}$ the function $f_{\omega }:\overline{G}\rightarrow \mathbb{C} $,%
\begin{equation*}
{\label{equation:fg}}f_{\omega }\left( \lambda \right) =\frac{1}{\omega
-f\left( \lambda \right) },\left( \forall \right) \lambda \in \overline{G},
\end{equation*}
is holomorphic, we can define $f_{\omega }\left( T\right) $, where 
\begin{equation}
{\label{equation:fg1}} f_{\omega }\left( T\right) =\frac{1}{2\pi i}%
\int_{\Gamma }\frac{1}{\omega -f\left( \lambda \right) }R(\lambda
,T)d\lambda =R(\omega ,f(T)).
\end{equation}

If we chose a Cauchy domain $G^{\prime }$ for the pair ($\sigma _{W}(Q_{\mathcal{P}%
},f(T)),D_{g})$ with the boundary $\Gamma ^{\prime }$ such that $\overline{%
f(G)}\subset G^{\prime }$, then $f(\Gamma )\cap \Gamma ^{\prime }=\varnothing $, so
we can define the function given by (\ref{equation:fg1}) for all $\lambda
\in \Gamma $ and $\omega \in \Gamma ^{\prime }$. Therefore, from relation (%
\ref{equation:fg1}) and Cauchy formula results%
\begin{equation*}
g(f\left( T\right) ) =\frac{1}{2\pi i}\int_{\Gamma ^{\prime }}g(\omega
)R(\omega ,f(T))d\omega =
\end{equation*}
\begin{equation*}
=\frac{1}{2\pi i}\int_{\Gamma ^{\prime }}g(\omega )\left( \frac{1}{2\pi i}%
\int_{\Gamma }\frac{1}{\omega -f\left( \lambda \right) }R(\lambda
,T)d\lambda \right) d\omega =
\end{equation*}
\begin{equation*}
=\frac{1}{2\pi i}\int_{\Gamma }R(\lambda ,T)\left( \frac{1}{2\pi i}%
\int_{\Gamma ^{\prime }}\frac{g(\omega )}{\omega -f\left( \lambda \right) }%
d\omega \right) d\lambda =\frac{1}{2\pi i}\int_{\Gamma }g(f(\lambda
))R(\lambda ,T)d\lambda =
\end{equation*}
\begin{equation*}
=\frac{1}{2\pi i}\int_{\Gamma }\left( g\circ f\right) (\lambda )R(\lambda
,T)d\lambda =\left( g\circ f\right) (T)
\end{equation*}
\end{proof}

\begin{lemma}
Assume that $\mathcal{P}\in C_{0}(X)$ and $T\in (Q_{\mathcal{P}}(X))_{0}$. If $f$ is a
holomorphic function on the open set $D\supset \sigma _{W}(Q_{\mathcal{P}%
},T) $ and $f\left( \lambda \right) =\sum\limits_{k=0}^{\infty }\alpha
_{k}\lambda ^{k}$ on $D$, then $f\left( T\right)
=\sum\limits_{k=0}^{\infty }\alpha _{k}T^{k}$.
\end{lemma}
\begin{proof}
For $\epsilon >0$ sufficiently small the power series $\sum\limits_{k=0}^{\infty }\alpha _{k}\lambda ^{k}$ converges uniformly on the boundary $\Gamma 
$ of the disc $D=\{\;\lambda \mid |\lambda |=|\sigma _{W}(Q_{\mathcal{P}%
},T)|+\epsilon \}$.

From corollary \ref{corollary:fc} results that%
\begin{equation*}
f(T)=\frac{1}{2\pi i}\int_{\Gamma }f\left( \lambda \right) R\left( \lambda
,T\right) d\lambda =\frac{1}{2\pi i}\int_{\Gamma }\left(
\sum\limits_{k=0}^{\infty }\alpha _{k}\lambda ^{k}\right) R\left( \lambda
,T\right) d\lambda =
\end{equation*}%
\begin{equation*}
=\frac{1}{2\pi i}\sum\limits_{k=0}^{\infty }\alpha _{k}\int_{\Gamma
}\lambda ^{k}R\left( \lambda ,T\right) \;d\lambda
\end{equation*}

Since for every $|\lambda |$ $>$ $|\sigma _{W}(Q_{\mathcal{P}},T)|$ we have $%
R\left( \lambda ,T\right) $=$\sum\limits_{k=0}^{\infty }\frac{T^{k}}{\lambda
^{k+1}}$, so from Cauchy formula results
\begin{center}
$f(T)=\frac{1}{2\pi i}\sum\limits_{k=0}^{\infty }\alpha _{k}\int_{\Gamma
}\lambda ^{k}\left( \sum\limits_{n=0}^{\infty }\frac{T^{n}}{\lambda ^{n+1}}%
\right) \;d\lambda =\sum\limits_{k=0}^{\infty }\alpha _{k}T^{k}.$
\end{center}
\end{proof}

\begin{corollary}
If $\mathcal{P}\in \mathcal{C}_{0}(X)$ and $T\in (Q_{\mathcal{P}}(X))_{0}$,
then $\exp T=\sum\limits_{k=0}^{\infty }\frac{{T}^{k}}{k!}$.
\end{corollary}

\begin{definition}
If $\mathcal{P}\in \mathcal{C}_{0}(X)$ and $T\in (Q_{\mathcal{P}}(X))_{0}$,
then a subset of $\sigma _{W}(Q_{\mathcal{P}},T)$ which is both open and
closed in $\sigma _{W}(Q_{\mathcal{P}},T)$ is called spectral set of $T$.
\end{definition}

Denote by $\delta _{T}$ the class of spectral sets of $T$.

\begin{proposition}
If $\mathcal{P}\in \mathcal{C}_{0}(X)$ and $T\in (Q_{\mathcal{P}}(X))_{0}$,
then for each spectral set $H\in \delta _{T}$ there exists a unique
idempotent $T_{H}\in Q_{\mathcal{P}}(X)$ with the following properties:
\begin{enumerate}
\item 
$T_{H}S=ST_{H}$, whenever $S\in Q_{\mathcal{P}}(X)$ and $ST=TS$;
\item 
$T_{\varnothing }$ is the null element of $Q_{\mathcal{P}}(X)$;
\item 
$T_{H\cap K}=T_{H}T_{K},\left( \forall \right) H,K\in \delta _{T}$;
\item 
$T_{H\cup K}=T_{H}+T_{K}$, for each $H,K\in \delta _{T}$ with the
property $H\cap K=\Phi $.
\end{enumerate}
\end{proposition}
\begin{proof}
First we make the observation that for each set $H\in \delta _{T}$\ there
exists an unique germs $\tilde{f}_{H}\in \mathcal{A}(\sigma _{W}(Q_{\mathcal{%
P}},T))$ with the property $(H)$, where 
\begin{equation*}
(H)\left\{ 
\begin{array}{l}
\text{for every pair }(D,D^{\prime })\text{ of relatively compact open
complex sets which  } \\ 
\text{ satiesfies the conditions } \\ 
\hspace{5.0188pc}H\subset D,~\sigma _{W}(Q_{\mathcal{P}},T)\subset D^{\prime
}\text{ and }D\cap D^{\prime }=\Phi \\ 
\text{ then there exists }f_{H}\in \tilde{f}_{H}\text{ such that }%
f_{H}/_{D}=1\text{ and }f_{H~}/_{D^{\prime }}=0.%
\end{array}%
\right.
\end{equation*}

If $\Gamma $ is Cauchy boundary for the pair $($ $\overline{H},D)$ (the
closure of$\ H$ is taken in the topology of $\mathbb{C}$) then by definition
consider that 
\begin{equation*}
T_{H}=F_{T}(\tilde{f}_{H})=\frac{1}{2\pi i}\int_{\Gamma }R(\lambda
,T)d\lambda .
\end{equation*}
1) If $ST=TS$, then $SR(\lambda ,T)=R(\lambda ,T)S$, so $T_{H}S=ST_{H}$.\\
2) Results from the definition of $T_{H}$.\\
3) Let $H,K\in $ $\delta _{T}$ and $\tilde{f}_{H},\tilde{f}_{K},\tilde{f}%
_{H\cap K}\in \mathcal{A}(\sigma _{W}(Q_{\mathcal{P}},T))$\ which verifies
the properties $(H)$, $(K)$, respectively $(H\cap K)$.

Assume that the pair $(D,D^{\prime })$ and $(G,G^{\prime })$ are like in $%
(H) $ and $(K)$ properties. Then there exists $f\in \tilde{f}_{H}$, such
that $f/_{D}=1$ and $f/_{D^{\prime }}=0$, and $g\in \tilde{f}_{K}$, such
that $g/_{G}=1$ and $g/_{G^{\prime }}=0$. It is obvious that $fg/_{D\cap
G}=1 $ and $fg/_{D^{\prime }\cap G^{\prime }}=0$, so $fg\in \tilde{f}_{H\cap
K}$ and%
\begin{equation*}
T_{H\cap M}=F_{T}(\tilde{f}_{H\cap K})=F_{T,D\cap G}(fg)=F_{T,D\cap
G}(f)F_{T,D\cap G}(g)=
\end{equation*}%
\begin{equation*}
=F_{T,D}(f)F_{T,G}(g)=F(\tilde{f}_{H})F(\tilde{f}_{K})=T_{H}T_{M}.
\end{equation*}

4) We consider the notations made above and with supplementary conditions
that $D\cap G=\varnothing $ and $D~^{\prime }\cap G~^{\prime }=\varnothing $, since $H\cap
M=\varnothing $. Then%
\begin{equation*}
f(\lambda )+g(\lambda )=\left\{ 
\begin{array}{l}
1,\text{\ if }\lambda \in D\cup G, \\ 
0,\text{ \ for }\lambda \in D~^{\prime }\cup G~^{\prime },%
\end{array}%
\right.
\end{equation*}

Therefore, if $ \tilde{f}_{H\cup K}\in \mathcal{A}(\sigma _{W}(Q_{\mathcal{P%
}},T))$ has the property $(H\cup K)$, then $\ f+g\in \tilde{f}_{H\cup K}$, so%
\begin{equation*}
T_{H\cup M}=F_{T}(\tilde{f}_{H\cup K})=F_{T,D\cup G}(f+g)=F_{T,D\cup
G}(f)+F_{T,D\cup G}(g)=
\end{equation*}%
\begin{equation*}
=F_{T,D}(f)+F_{T,G}(g)=F_{T}(\tilde{f}_{H})+F_{T}(\tilde{f}_{K})=T_{H}+T_{M}.
\end{equation*}
\end{proof}

\begin{corollary}
If $\mathcal{P}\in \mathcal{C}_{0}(X)$ and $T\in (Q_{\mathcal{P}}(X))_{0}$,
then for every pair of spectral sets spectral set $H,K\in \delta _{T}$,
which have the properties $H\cap K=\varnothing $ and $H\cup K=\sigma _{W}(Q_{%
\mathcal{P}},T)$ , we have 
\begin{equation*}
T_{H}+T_{K}=I\text{ and }T_{H}T_{K}=O.
\end{equation*}
\end{corollary}

\begin{remark}
From proposition \ref{proposition:reg1} results that $T_{H}\in (Q_{\mathcal{P%
}}(X))_{0}$, for each $H\in \delta _{T}$
\end{remark}

\begin{lemma}
Assume that $\mathcal{P}\in \mathcal{C}_{0}(X)$ and $T\in (Q_{\mathcal{P}%
}(X))_{0}$. If $F\subset \mathbb{C}$ has the property $dist(\sigma _{W}(Q_{%
\mathcal{P}},T),F)>\varepsilon_{0}>0$, then for each $p\in \mathcal{P}$ there
exists $r_{p}>0$ such that 
\begin{equation*}
\hat{p}(R(\lambda ,T)^{n})\leq \frac{r_{p}}{\varepsilon_{0} ^{n}},~\left( \forall
\right) \lambda \in F,\left( \forall \right) n\in \mathbb{N}
\end{equation*}
\end{lemma}
\begin{proof}
Let be $\varepsilon\in(o,\varepsilon_{0})$, arbitrary fixed. If $D=\mathbb{C}\backslash \overline{F}$, then for the pair $(\sigma _{W}(Q_{\mathcal{P}},T)),D)$ there exists a Cauchy domain $G$ such that 
\begin{equation*}
|\lambda -\omega |>{\varepsilon_{0}-\varepsilon },~\left( \forall \right) \lambda
\in F,\left( \forall \right) \omega \in \overline{G}
\end{equation*}

If $\Gamma $ is boundary of $G$, then 
\begin{equation*}
\hat{p}(R(\lambda ,T)^{n})=\hat{p}\left( \frac{1}{2\pi i}\int_{\Gamma }\frac{%
R(\omega ,T)}{\left( \omega -\lambda \right) ^{n}}d\omega \right) \leq \frac{%
L(\Gamma )}{2\pi }\sup_{\omega \in \Gamma }\frac{\hat{p}(R(\omega ,T)^{n})}{%
\left\vert \omega -\lambda \right\vert ^{n}}<
\end{equation*}%
\begin{equation*}
<\frac{\frac{L(\Gamma )}{2\pi }\sup_{\omega \in \Gamma }\hat{p}(R(\omega
,T)^{n})}{(\varepsilon_{0}-\varepsilon ) ^{n}}.
\end{equation*}
Since $\varepsilon$ is arbitrary results that for $r_{p}=\frac{L(\Gamma )}{2\pi }\sup_{\omega \in \Gamma }\hat{p}%
(R(\omega ,T)^{n})$ the lemma is proved.
\end{proof}

\begin{theorem}
Let $\mathcal{P}\in \mathcal{C}_{0}(X)$ and $T\in (Q_{\mathcal{P}}(X))_{0}$.
If $D$ is an open relatively compact set which contains the set $\sigma
_{W}(Q_{\mathcal{P}},T),f\in \mathcal{O}(D)$ and $S\in (Q_{\mathcal{\ P}%
}(X))_{0}$, such that $r_{\mathcal{P}}(S)<dist(\sigma _{W}(Q_{\mathcal{P}%
},T),\mathbb{C}\backslash D)$ and $TS=ST$, then the following statements ar
true:
\begin{enumerate}
\item 
$\sigma _{W}(Q_{\mathcal{P}},T+S)\subset D$;
\item 
$f(T+S)=\sum\limits_{n\geq 0}\frac{f^{\left( n\right) }\left( T\right) 
}{n!}S^{n}$.
\end{enumerate}
\end{theorem}
\begin{proof}
Let $d,d_{1}>0$ such that 
\begin{equation*}
r_{\mathcal{P}}(S)<d_{1}<d<dist\left( \sigma _{W}(Q_{\mathcal{P}},T),\mathbb{%
C}\backslash D\right) .
\end{equation*}

If $\Gamma _{1}=\{\lambda \in \mathbb{C~}|~\left\vert \lambda \right\vert
=d_{1}\}$, then for each $p\in \mathcal{P}$ and each $n\in \mathbb{N}$ we
have 
\begin{equation}
{\label{equation:fg 2}} \hat{p}(S^{n})=\hat{p}\left( \frac{1}{2\pi i}%
\int_{\Gamma }\lambda ^{n}R(\lambda ,S)d\lambda \right) \leq \frac{L(\Gamma
_{1})}{2\pi }\sup_{\omega \in \Gamma _{1}}\left( \left\vert \lambda
^{n}\right\vert \hat{p}(R(\lambda ,S)\right) \leq
\end{equation}%
\begin{equation*}
\leq \frac{L(\Gamma _{1})}{2\pi }\sup_{\omega \in \Gamma _{1}}\hat{p}%
(R(\lambda ,S)\sup_{\omega \in \Gamma _{1}}\left\vert \lambda \right\vert
^{n}\leq k_{p}d_{1}^{n}
\end{equation*}
where $k_{p}=\frac{L(\Gamma _{1})}{2\pi }\sup_{\omega \in \Gamma _{1}}\hat{p}%
(R(\lambda ,T)$.

Moreover, the previous lemma implies that for each $p\in \mathcal{P}$ there
exists $r_{p}>0$ such that%
\begin{equation}{\label{equation:fg 3}}
\hat{p}(R(\lambda ,T)^{n+1})\leq\frac{r_{p}}{d^{n+1}}%
,~\left( \forall \right) \lambda \in \mathbb{C}\backslash D,\left( \forall
\right) n\in \mathbb{N}
\end{equation}
so from relation (\ref{equation:fg 2}) and (\ref{equation:fg 3}) results
that 
\begin{equation}{\label{equation:fg 4}}
\hat{p}(R(\lambda ,T)^{n+1}S^{n})=\hat{p}(R(\lambda
,T)^{n+1})\hat{p}(S^{n})\leq \frac{k_{p}r_{p}}{d_{1}}\left( \frac{d_{1}}{d}%
\right) ^{n+1}
\end{equation}
for every $p\in \mathcal{P}$, $n\in \mathbb{N}$ and $\lambda \in \mathbb{C}%
\backslash D$. Since $\frac{d_{1}}{d}<1$ the relation (\ref{equation:fg 4})
prove that the series $\sum_{n=1}^{\infty }R(\lambda ,T)^{n+1}S^{n}$
converge uniformly on $\mathbb{C}\backslash D$.

From equalities
\begin{equation*}
(\lambda I-T-S)\sum_{n=1}^{\infty }R(\lambda
,T)^{n+1}S^{n}=\sum_{n=1}^{\infty }R(\lambda ,T)^{n+1}S^{n}(\lambda I-T-S)=
\end{equation*}%
\begin{equation*}
=\sum_{n=1}^{\infty }R(\lambda ,T)^{n}S^{n}-\sum_{n=1}^{\infty }R(\lambda
,T)^{n+1}S^{n+1}=I
\end{equation*}
results that $\lambda I-T-S$ is invertible in $Q_{\mathcal{P}}(X)$, for all $%
\lambda \in \mathbb{C}\backslash D$, and 
\begin{equation}{\label{equation:fg 5}} 
R(\lambda ,T+S)=\sum_{n=1}^{\infty }R(\lambda,T)^{n+1}S^{n}.
\end{equation}

Therefore the relation (\ref{equation:fg 4}) implies that the set $%
\{R(\lambda ,T+S)|\lambda \in \mathbb{C}\backslash D\}$ is bounded in $Q_{%
\mathcal{P}}(X)$, so $\sigma _{W}(Q_{\mathcal{P}},T+S)\subset D$.

If $\Gamma $ is a Cauchy boundary for the pair $(\sigma _{W}(Q_{\mathcal{P}%
},T+S),D)$, then from (\ref{equation:fg 5}) and lemma \ref{lemma:qb1}
results 
\begin{equation*}
{\label{equation:fg 6}}f(\lambda ,T+S)=\frac{1}{2\pi i}\int_{\Gamma
}f(\lambda )R(\lambda ,T+S)d\lambda =
\end{equation*}%
\begin{equation*}
=\sum_{n=1}^{\infty }\left( \frac{1}{2\pi i}\int_{\Gamma }f(\lambda
)R(\lambda ,T)^{n+1}d\lambda \right) S^{n}=
\end{equation*}%
\begin{equation*}
=\sum_{n=1}^{\infty }\left( \frac{1}{2\pi i}\frac{(-1)^{n}}{n!}\int_{\Gamma
}f(\lambda )\frac{d^{n}}{d \lambda ^{n}} R(\lambda ,T)d\lambda \right) S^{n}=
\end{equation*}%
\begin{equation*}
=\sum_{n=1}^{\infty }\left( \frac{1}{2\pi i}\frac{1}{n!}\int_{\Gamma
}f^{(n)}(\lambda )R(\lambda ,T)d\lambda \right) S^{n}=\sum_{n=1}^{\infty }%
\frac{f^{(n)}(T)}{n!}S^{n}
\end{equation*}
\end{proof}

\begin{corollary}
Let $\mathcal{P}\in \mathcal{C}_{0}(X)$ and $T\in (Q_{\mathcal{P}}(X))_{0}$.
If $S\in Q_{\mathcal{P}}(X)$ is $\mathcal{P}$-quasnilpotent, such that $TS=ST
$, then  
\begin{equation*}
\tilde{f}(T+S)=\sum\limits_{n\geq 0}\frac{\tilde{f}^{\left( n\right) }\left(
T\right) }{n!}S^{n},~(\forall)\tilde{f}\in \mathcal{A}(\sigma _{W}(Q_{\mathcal{P}},T))
\end{equation*}
\end{corollary}

\section{Locally bounded operators}

\begin{theorem}
If $T$ $\in \mathcal{LB}(X)$, then 
\begin{equation*}
\sigma (T)=\cap \{\sigma (B_{\mathcal{P}},T){|}T\in B_{\mathcal{P}}(X)\}
\end{equation*}
\end{theorem}
\begin{proof}
Assume that $\lambda \notin $ $\sigma (T)$. From \ lemma \ref{lemma:k2}
results that there exists some calibration $\mathcal{P}$ $\in \mathcal{C}(X)$
such that $(\lambda I-T)^{-1},T\in B_{\mathcal{P}}(X)$, i.e. $\lambda \in
\rho (B_{\mathcal{P}},T)$. From this observation results that $\lambda
\notin \cap \{\sigma (B_{\mathcal{P}},T){|}T\in B_{\mathcal{P}}(X)\}$, so 
\begin{equation*}
\cap \{\sigma (B_{\mathcal{P}},T){|}T\in B_{\mathcal{P}}(X)\}\subset \sigma
(T)
\end{equation*}

The reverse inclusion is obvious.
\end{proof}

\begin{theorem}
If $T\in \mathcal{LB}(X)$, then there exists $\mathcal{P}\in \mathcal{C}(X)$
such that $T\in B_{\mathcal{P}}(X)$ and 
\begin{equation*}
\left\vert \sigma \left( T\right) \right\vert ={|}\sigma (B_{\newline
\mathcal{P}},T){|}=\lim\limits_{n\rightarrow \infty }({|}{|}T^{n}{|}{|}_{%
\mathcal{P}})^{1/n}=r_{lb}(T)
\end{equation*}
\end{theorem}
\begin{proof}
Let $\mathcal{P}$ $\in \mathcal{C}_{0}(X)$ such that $T\in B_{\mathcal{P}%
}(X) $. For every $p\in\mathcal{P}$ and $\epsilon>0$ we consider
\begin{equation*}
S_{p}(0_{X},\epsilon )=\left \{ x\in X \vert p(x)<\epsilon \right \}
\end{equation*}

If $\nu >r_{lb}(T)$, then there exists $\mu \in (r_{lb}(T),\nu )$ such that
the sequence $\left\{\frac{T^{n}}{\mu ^{n}}\right\}_{n}$ converges to zero on a zero neighborhood. Then, there exists a zero neighborhood $U$ with the property that for each
neighborhood $V$ there exists an index $n_{V}\in \mathbf{N}$ such that 
\begin{equation*}
\frac{T^{n}}{\mu ^{n}}(U)\subset V,\left( \forall \right) n\geq n_{V}
\end{equation*}

Since $\mathcal{P}\in \mathcal{C}_{0}(X)$ we can assume without lost of the
generality of the proof that there exists $p_{0}\in \mathcal{P}$ and $%
\epsilon >0$ such that%
\begin{equation*}
\frac{T^{n}}{\mu ^{n}}(S_{p_{0}}\left( 0_{X},\epsilon \right) )\subset
V,\left( \forall \right) n\geq n_{V}.
\end{equation*}

Moreover, since $T\in \mathcal{LB}(X)$ there exists some seminorm $p_{1}\in 
\mathcal{P}$ such that the set $TS_{p_{1}}\left( 0_{X },\epsilon \right) $
is bounded, so for each zero neighborhood $V$ there exists $\alpha _{V}>0$
such that 
\begin{equation}
{\label{equation:lb1}}TS_{p_{1}}\left( 0_{X},\epsilon \right) \subset \alpha
_{V}V.
\end{equation}

Let be $n_{1}$ the index for which we have%
\begin{equation*}
\frac{T^{n}}{\mu ^{n}}S_{p_{0}}\left( 0_{X },\epsilon \right) \subset
S_{p_{1}}\left( 0_{X},\epsilon \right) ,\left( \forall \right) n\geq n_{1}.
\end{equation*}

Then for each zero neighborhood $V$ there exists $\alpha _{V}>0$ such that 
\begin{equation*}
\frac{T^{n+1}}{\mu ^{n+1}}S_{p_{0}}\left( 0_{X },\epsilon \right) =\frac{T}{%
\mu }\left( \frac{T^{n}}{\mu ^{n}}S_{p_{0}}\left( 0_{X},\epsilon \right)
\right) \subset \frac{1}{\mu }TS_{p_{1}}\left( 0_{X},\epsilon \right)
\subset \frac{\alpha _{V}}{\mu }V,\left( \forall \right) n\geq n_{1}.
\end{equation*}
so if $\beta _{V}=\mu ^{-1}\alpha _{V}$ and $n_{2}=n_{1}+1$, then%
\begin{equation*}
\frac{T^{n}}{\mu ^{n}}S_{p_{0}}\left( 0_{X},\epsilon \right) \subset \beta
_{V}V,\left( \forall \right) n\geq n_{2}.
\end{equation*}

Assume that $V=S_{q}\left( 0_{X},1\right) $, where $q\in \mathcal{P}$.
Therefore, for each $q\in \mathcal{P}$ we can find $\beta _{q}>0$ such that 
\begin{equation*}
\frac{T^{n}}{\mu ^{n}}S_{p_{0}}\left( 0_{X},\epsilon \right) \subset \beta
_{q}S_{q}\left( 0_{X},1\right) ,\left( \forall \right) n\geq n_{2}.
\end{equation*}

This implies that if $p_{0}\left( x\right) <1$ then 
\begin{equation*}
q\left( \frac{\epsilon }{\beta _{q}}\frac{T^{n}}{\mu ^{n}}x\right) <1,\left(
\forall \right) n\geq n_{2},
\end{equation*}
so 
\begin{equation*}
q\left( \frac{T^{n}}{\mu ^{n}}x\right) \leq \delta _{q}p_{0}\left( x\right)
,\left( \forall \right) n\geq n_{2},\left( \forall \right) x\in X,
\end{equation*}
where $\delta _{q}=\epsilon ^{-1}\beta _{q}$.

Using relation (\ref{equation:lb1}) by the same method we can find for each $%
q\in \mathcal{P}$ a pozitiv number $\tau _{q}$ such that 
\begin{equation*}
q\left( \text{T}x\right) \leq \tau _{q}p_{1}\left( x\right) ,\left( \forall
\right) x\in X.
\end{equation*}

Since $\mathcal{P}\in \mathcal{C}_{0}(X)$ there exists $p_{2}\in \mathcal{%
P}$ such that $p_{0}\leq p_{2}$ and $p_{1}\leq p_{2}$. Therefore, for each $%
q\in \mathcal{P}$ we have%
\begin{equation*}
q\left( \frac{T^{n}}{\mu ^{n}}x\right) \leq \gamma _{q}p_{2}\left( x\right)
,\left( \forall \right) n\geq n_{2},\left( \forall \right) x\in X,
\end{equation*}%
\begin{equation*}
q\left( Tx\right) \leq \gamma _{q}p_{2}\left( x\right) ,\left( \forall
\right) x\in X.
\end{equation*}
where $\gamma _{q}=\max \left\{ \tau _{q},\delta _{q}\right\} $.

For each $q\in \mathcal{P}$ consider the application $q^{\prime}:X\rightarrow \mathbb{R}$ given by relation
\begin{equation*}
q^{\prime }(x)=\max \left\{ q\left( x\right),\gamma _{q}p_{2}\left(
x\right) \right\}, \left( \forall
\right) x\in X.
\end{equation*}%
It is easy to observe that $q^{\prime }$ is a seminorm on $X$ and $\mathcal{P%
}^{\prime }=\{q~^{\prime }|q\in \mathcal{P}\}\in \mathcal{C}(X)$.

If $c_{0}=\max \left\{ 1,\gamma _{p_{2}}\right\} $, then for each $q^{\prime
}\in \mathcal{P}^{\prime} $ we have%
\begin{equation*}
q^{\prime }(Tx)=\max \left\{ q\left( Tx\right) ,\gamma _{q}p_{2}\left(
Tx\right) \right\} \leq \max \left\{ \gamma _{q}p_{2}\left( x\right) ,\gamma
_{q}\gamma _{p_{2}}p_{2}\left( x\right) \right\} =
\end{equation*}%
\begin{equation*}
=\gamma _{q}p_{2}\left( x\right) \max \left\{ 1,\gamma _{p_{2}}\right\} \leq
c_{0}q^{\prime }\left( x\right) ,\left( \forall \right) x\in X.
\end{equation*}%
\begin{equation*}
q^{\prime }\left( \frac{T^{n}}{\mu ^{n}}x\right) =\max \left\{ q\left( \frac{%
T^{n}}{\mu ^{n}}x\right) ,\gamma _{q}p_{2}\left( \frac{T^{n}}{\mu ^{n}}%
x\right) \right\} \leq
\end{equation*}%
\begin{equation*}
\leq \max \left\{ \gamma _{q}p_{2}\left( x\right) ,\gamma _{q}\gamma
_{p_{2}}p_{2}\left( x\right) \right\} =\gamma _{q}p_{2}\left( x\right) \max
\left\{ 1,\gamma _{p_{2}}\right\} \leq c_{0}q^{\prime }\left( x\right)
\end{equation*}
for every $n\geq n_{2}$ and every $x\in X$, so $T\in B_{\mathcal{P}^{\prime}}(X)$.

From inequality%
\begin{equation*}
q^{\prime }\left( \frac{T^{n}}{\mu ^{n}}x\right) \leq c_{0}q^{\prime }\left(
x\right) ,\left( \forall \right) n\geq n_{2},\left( \forall \right) x\in
X,\left( \forall \right) q^{\prime }\in \mathcal{P}^{\prime} ,
\end{equation*}
results that%
\begin{equation*}
\left\Vert \frac{T^{n}}{\mu ^{n}}\right\Vert _{\mathcal{P}^{\prime}}\leq
c_{0},\left( \forall \right) n\geq n_{2},
\end{equation*}
or equivalently%
\begin{equation*}
\left\Vert T^{n}\right\Vert _{\mathcal{P}^{\prime}}\leq c_{0}\mu ^{n},\left( \forall
\right) n\geq n_{2}.
\end{equation*}

Therefore%
\begin{equation*}
\underset{n\rightarrow \infty }{\limsup }({|}{|}T^{n}{|}{|}_{\mathcal{P}^{\prime}%
}^{ 1/n}\leq \mu <\nu .
\end{equation*}%
and since $\nu >r_{lb}(T)$ is arbitrary chosen results%
\begin{equation*}
\underset{n\rightarrow \infty }{\limsup }({|}{|}T^{n}{|}{|}_{\mathcal{P}%
^{\prime }})^{1/n}\leq r_{lb}(T).
\end{equation*}

From inclusion $\sigma (T)\subset \sigma (B_{\mathcal{P}^{\prime }},T)$,
propositions \ref{proposition:raza spectrala lb} and \ref%
{proposition:univ.bound.1} results that 
\begin{equation*}
r_{lb}(T)=\left\vert \sigma \left( T\right) \right\vert \leq |\sigma (B_{%
\mathcal{P}^{\prime}},T)|\leq \underset{n\rightarrow \infty }{\liminf }({|}{|}T^{n}{|%
}{|}_{\mathcal{P}^{\prime}})^{1/n},
\end{equation*}%
so%
\begin{equation*}
\left\vert \sigma \left( T\right) \right\vert ={|}\sigma (B_{\mathcal{P}^{\prime}%
},T)|=\lim\limits_{n\rightarrow \infty }({|}{|}T^{n}{|}{|}_{{\mathcal{P}}^{\prime}})^{1/n}=r_{lb}(T).
\end{equation*}
\end{proof}

For a locally bounded operator $T$ on a locally convex space $X$ we can
define the main subsets of the spectrum $\sigma _{lb}(T)$: the point
spectrum $\sigma _{p}(T)$, the continuous spectrum $\sigma _{c}(T)$, the
residual spectrum $\sigma _{r}(T))$, respectively the approximate spectrum $%
\sigma _{a}(T)$.

\begin{definition}
If $X$ is a locally convex space and $T\in \mathcal{LB}(X)$, then:
\begin{enumerate}
\item 
$\lambda \in \sigma _{p}(T)$ if and only if $\ker (\lambda I-T)\neq
\left\{ 0_{X}\right\} $;
\item 
$\lambda \in \sigma _{c}(T)$ if and only if $(\lambda
I-T)^{-1}$ exists on the set $\text{Im}\lambda (I-T)$ which is dense in X and 
$(\lambda I-T)X\neq X$;
\item 
$\lambda \in \sigma _{r}(T)$ if and only if  $(\lambda
I-T)^{-1}$ exists on the set $\text{Im}(\lambda I-T)$, which is not dense in 
$X$;
\item 
$\lambda \in \sigma _{a}(T)$ if and only if for each $\mathcal{P}\in 
\mathcal{C}(X)$, such that $T\in B_{\mathcal{P}}(X)$, and for every $c>0$
there exists a seminorm $p\in \mathcal{P}$ and an element $x\in X$ such that 
$\ p\left( \left( \lambda I-T\right) x\right) <cp\left( x\right) $.
\end{enumerate}
\end{definition}

\begin{remark}
From previous definition results that the sets $\sigma _{p}(T),\sigma _{c}(T)$
and $\sigma _{r}(T)$ are disjoint and 
\begin{equation*}
\sigma (T)=\sigma _{p}(T)\cup \sigma _{c}(T)\cup \sigma _{r}(T).
\end{equation*}
\end{remark}

\begin{theorem}
If $T$ $\in \mathcal{LB}(X)$, then
\begin{enumerate}
\item 
$\sigma _{r}(T)\cup \sigma _{a}(T)=\sigma (T)$;
\item 
$\sigma _{p}(T)\subset \sigma _{a}(T)$ and $\sigma _{c}(T)\subset \sigma
_{a}(T)$.
\end{enumerate}
\end{theorem}
\begin{proof}
1) Denote by $\sigma _{r}(T)^{c}$ and $\sigma _{a}(T)^{c}$ the
complements of the sets $\sigma _{r}(T)$, respectively $\sigma _{a}(T)$. Let $%
\lambda \in \sigma _{r}(T)^{c}\cap \sigma _{a}(T)^{c}$ and $y\in X$.

From the condition $\lambda \notin \sigma _{a}(T)$ results that there exists 
$\mathcal{P}\in \mathcal{C}(X)$, such that $T\in B_{\mathcal{P}}(X)$, and $%
c_{0}>0$ such that 
\begin{equation*}
p\left( \left( \lambda I-T\right) x\right) \geq c_{0}p\left( x\right)
,\left( \forall \right) x\in X,\left( \forall \right) p\in \mathcal{P},
\end{equation*}%
so $\ker (\lambda I-T)=\left\{ 0_{X}\right\} $. Therefore the operator $%
(\lambda I-T)^{-1}$ exists on the set $\text{Im}(\lambda I-T)$, which is
dense in X. The previous inequality is equivalent with 
\begin{equation}
{\label{equation:lb2}} p\left( \left( \lambda I-T\right) ^{-1}z\right) \leq
c_{0}^{-1}p\left( z\right) ,\left( \forall \right) z\in \text{Im}(\lambda
I-T),\left( \forall \right) p\in \mathcal{P}.
\end{equation}

Since $\text{Im}(\lambda I-T)$ is dense in X, then for each $y\in X$ there exists an sequence $%
\left\{ x_{\delta }\right\} _{\delta }$ such that $y_{\delta }=\lambda
x_{\delta }-Tx_{\delta }$ converges to $y$.

From previous observation results that the operator $\left( \lambda
I-T\right) ^{-1}$ exists on the set $\text{Im}(\lambda I-T)$ and is
continuous in the sense of the relation (\ref{equation:lb2}). Therefore the
sequence $x_{\delta }=(\lambda I-T)^{-1}y_{\delta }$ converges in $X$ to an
unique element $x$ and from continuity of $(\lambda I-T)^{-1}$ results that $%
(\lambda I-T)x=y$

Therefore, $(\lambda I-T)X=X$, and from (\ref{equation:lb2}) results that $%
\left( \lambda I-T\right) ^{-1}\in B_{\mathcal{P}}(X)$, i.e. $\lambda
\notin \sigma (T)$. This implies that $\sigma (T)\subset \sigma _{r}(T)\cup
\sigma _{a}(T)$.

For reverse inclusion we must prove that $\sigma _{a}(T)\subset \sigma (T)$.
If $\lambda \in \sigma _{a}(T)$, then for each $\mathcal{P}\in \mathcal{C}(X)$, with the property $T\in B_{\mathcal{P}}(X)$, and every $c>0$ there exists
a seminorm $p\in \mathcal{P}$ and an element $x\in X$ such that%
\begin{equation*}
p\left( \left( \lambda I-T\right) x\right) <cp\left( x\right)
\end{equation*}

If we assume that $\lambda \notin \sigma (T)$, then there exists a calibration $\mathcal{P}%
\in \mathcal{C}(X)$ such that $(\lambda I-T)^{-1},T\in B_{\mathcal{P}}(X)$.
From lemmas \ref{lemma:op.semin.} and \ref{lemma:1} results that 
\begin{equation*}
p\left( \left( \lambda I-T\right) ^{-1}x\right) \leq {|}{|}\left( \lambda
I-T\right) ^{-1}{|}{|}_{\mathcal{P}}p\left( x\right) ,\left( \forall \right)
x\in X,\left( \forall \right) p\in \mathcal{P},
\end{equation*}
or equivalently that for $y=\left( \lambda I-T\right) ^{-1}x$ we have%
\begin{equation*}
p\left(y\right) \leq {|}{|}\left( \lambda
I-T\right) ^{-1}{|}{|}_{\mathcal{P}}p\left( \left( \lambda
I-T\right)y\right) ,\left( \forall \right)
x\in X,\left( \forall \right) p\in \mathcal{P},
\end{equation*}
which contradict the supposition we made. \newline
2) Since the sets $\sigma _{p}(T) $, $\sigma _{c}(T)$ and $\sigma _{r}(T)$
are disjoint two of each, the property follows from (1)
\end{proof}

\begin{theorem}
If $T\in \mathcal{LB}(X)$, then

\begin{enumerate}
\item 
the boundary $\partial \sigma (T)$ of the spectral set $\sigma (T)$ is
included in $\sigma _{a}(T)$;
\item 
the set $\sigma _{a}(T)$\ is bounded and closed.
\end{enumerate}
\end{theorem}
\begin{proof}
1) Let $\lambda _{0}\in \partial \sigma (T)$ and $\epsilon >0$. Then there
exists $\lambda \in \rho (T)=\rho _{lm}(T)$ such that $\mid \lambda -\lambda
_{0}\mid <\frac{\epsilon }{2}$. From lemma \ref{lemma:k2} results that there
exists $\mathcal{P}\in \mathcal{C}(X)$ such that $(\lambda I-T)^{-1},T\in B_{%
\mathcal{P}}(X)$. Then from corollary \ref{corollary:1} we have
\begin{equation*}
\left\Vert \left( \lambda I-T\right) ^{-1}\right\Vert _{\mathcal{P}}\geq 
\frac{1}{d\left( \lambda \right) }>\frac{2}{\epsilon },
\end{equation*}
where $d\left( \lambda \right) $ is the distance from $\lambda $ to $\sigma
(T)$. Since 
\begin{equation*}
\left\Vert \left( \lambda I-T\right) ^{-1}\right\Vert _{\mathcal{P}}=\inf
\{M>0\left\vert p\left( \left( \lambda I-T\right) ^{-1}x\right) \leq \right.
Mp\left( x\right) ,\left( \forall \right) x\in X,\left( \forall \right) p\in 
\mathcal{P}\},
\end{equation*}
results that there exists $p_{0}\in \mathcal{P}$ and $x_{0}\in X$ such that 
\begin{equation*}
p_{0}\left( \left( \lambda I-T\right) ^{-1}x_{0}\right) >2\epsilon
^{-1}p_{0}\left( x_{0}\right) .
\end{equation*}

For $y_{0}=\left( \lambda I-T\right) ^{-1}x_{0}$ we have 
\begin{equation*}
p_{0}\left( \left( \lambda I-T\right) y_{0}\right) <2^{-1}\epsilon
p_{0}\left( y_{0}\right) .
\end{equation*}

If $\mathcal{P}\in \mathcal{C}(X)$ such that $T\in B_{\mathcal{P}}(X)$ and $%
(\lambda I-T)^{-1}\notin B_{\mathcal{P}}(X)$, then from the definition of
universally bounded operators there exists $p_{1}\in \mathcal{P}$ and $%
y_{1}\in X$ such that 
\begin{equation*}
p_{1}\left( \left( \lambda I-T\right) ^{-1}y_{1}\right) >2\epsilon
^{-1}p_{1}\left( y_{1}\right)
\end{equation*}

Therefore, for $y_{2}=\left( \lambda I-T\right) ^{-1}y_{1}$ we have 
\begin{equation*}
p_{1}\left( \left( \lambda I-T\right) y_{2}\right) <2^{-1}\epsilon
p_{1}\left( y_{2}\right)
\end{equation*}

In conclusion, for every $\mathcal{P}\in \mathcal{C}(X)$, such that $T\in B_{%
\mathcal{P}}(X)$, and for every $\epsilon >0$ there exists $q_{0}\in 
\mathcal{P}$ and $z_{0}\in X$ such that%
\begin{equation*}
q_{0}\left( \left( \lambda _{0}I-T\right) z_{0}\right) \leq q_{0}\left(
\left( \lambda _{0}I-T\right) z_{0}-\left( \lambda I-T\right) z_{0}\right)
+q_{0}\left( \left( \lambda I-T\right) z_{0}\right) =
\end{equation*}%
\begin{equation*}
=q_{0}\left( \left( \lambda _{0}-\lambda \right) z_{0}\right) +q_{0}\left(
\left( \lambda I-T\right) z_{0}\right) <\left( \frac{\epsilon }{2}+\frac{%
\epsilon }{2}\right) q_{0}\left( z_{0}\right) =\epsilon q_{0}\left(
z_{0}\right)
\end{equation*}

Therefore, $\lambda _{0}\in \sigma _{a}(T)$ and $\partial \sigma (T) \subset
\sigma _{a}(T)$.\newline
2)	Since a compact set has a nonempty boundary from (1) results that the set $%
\sigma _{a}(T)$ is nonmpty set.

If $\lambda _{0}\in\sigma _{a}(T)^{c}$, then there exists $\mathcal{P}\in 
\mathcal{C}(X)$, such that $T\in B_{\mathcal{P}}(X)$, and $c_{0}>0$ with the
property%
\begin{equation*}
p\left( \left( \lambda _{0}I-T\right) x\right) \geq c_{0}p\left( x\right)
,\left( \forall \right) x\in X,\left( \forall \right) p\in \mathcal{P},
\end{equation*}

For $\lambda \in \mathbb{C}$ with the property $\mid \lambda -\lambda
_{0}\mid <c_{0}/2$ we have 
\begin{equation*}
c_{0}p\left( x\right) \leq p\left( \left( \lambda _{0}I-T\right) x\right)
\leq
\end{equation*}%
\begin{equation*}
\leq p\left( \left( \lambda _{0}I-T\right) x+\left( \lambda I-T\right)
x\right) +p\left( \left( \lambda I-T\right) x\right) \leq
\end{equation*}%
\begin{equation*}
\leq p\left( \left( \lambda _{0}-\lambda \right) x\right) +p\left( \left(
\lambda I-T\right) x\right) =2^{-1}c_{0}p\left( x\right) +p\left( \left(
\lambda I-T\right) x\right)
\end{equation*}
for all $x\in X$ and $p\in \mathcal{P}$, so 
\begin{equation*}
2^{-1}c_{0}p\left( x\right) \leq p\left( \left( \lambda I-T\right) x\right)
,\left( \forall \right) x\in X,\left( \forall \right) p\in \mathcal{P}.
\end{equation*}

Therefore, $\lambda \in \sigma _{a}(T)$ $^{c}$ and 
\begin{equation*}
\left\{ \lambda \mid \mid \lambda -\lambda _{0}\mid <\frac{c_{0}}{2}\right\}
\subset \sigma _{a}(T)^{c}.
\end{equation*}

Since $\lambda _{0}\in \sigma _{a}(T)^{c}$ is arbitrary chosen results that
the set $\sigma _{a}(T)^{c}$ is open, so (2) is proved.
\end{proof}

\begin{theorem}
{\label{theorem:lb1}}If $T\in \mathcal{LB}(X)$ and $\mathcal{P}\in \mathcal{C%
}(X)$, such that $T\in B_{\mathcal{P}}(X)$ and $\sigma (Q_{\mathcal{P}},T)$
is closed, then $\sigma (Q_{\mathcal{P}},T)=\sigma _{W}(Q_{\mathcal{P}},T)$.
\end{theorem}
\begin{proof}
Let $\mathcal{P}=\left( p_{\alpha }\right) _{\alpha \in \Lambda }\in 
\mathcal{C}(X)$ such that $T\in B_{\mathcal{P}}(X)$. From definition of
Waelbroeck spectrum results that $\rho _{W}(Q_{\mathcal{P}},T)\subset \rho
(Q_{\mathcal{P}},T)$.

We prove that $\rho (Q_{\mathcal{P}},T)\subset \rho _{W}(Q_{\mathcal{P}},T)$. Let $\lambda _{0}\in \rho (Q_{\mathcal{P}},T)$ ( $\lambda _{0}\neq \infty $%
). Since $\sigma (Q_{\mathcal{P}},T)$ is closed set results that $\rho (Q_{%
\mathcal{P}},T)$ is open, so there exists $\varepsilon >0$ such that $%
D\left( \lambda _{0},\epsilon \right) \subset \rho (Q_{\mathcal{P}},T)$,
i.e. for every $\lambda \in D\left( \lambda _{0},\epsilon \right) $ the
operator $\lambda I-T$ is invertible and $(\lambda I-T)^{-1}\in Q_{\mathcal{P%
}}(X)$.

We will prove that there exists $\varepsilon _{0}>0$ such that $%
\{R(\lambda ,T){|}\lambda \in D\left( \lambda _{0},\epsilon _{0}\right) \}$
is a bounded set in $Q_{\mathcal{P}}(X)$.

First we study if $\sigma (Q_{\mathcal{P}},R(\lambda _{0},T))$ is bounded.
If $\mid \mu \mid >\epsilon ^{-1}$, then $\mid \mu \mid ^{-1}<\epsilon $ and 
$\lambda _{0}-\mu ^{-1}\in D\left( \lambda _{0},\epsilon \right) $. Therefore, results that the operator $(\lambda _{0}-\mu ^{-1})I-T$ is invertible and 
$((\lambda _{0}-\mu ^{-1})I-T)^{-1}\in Q_{\mathcal{P}}(X)$.

We need to prove the equality 
\begin{equation}{\label{equation:lb3}}
R(\mu ,R\left( \lambda _{0},T \right) )=\mu
^{-1}R\left( \lambda _{0}-\mu ^{-1},T \right) \left( \lambda _{0} I-T\right)
\end{equation}

Indeed%
\begin{equation*}
\mu ^{-1}R\left( \lambda _{0}-\mu ^{-1},T\right)
\left( \lambda _{0}I-T\right) (\mu I-R(\lambda _{0},T))=
\end{equation*}%
\begin{equation*}
=R\left( \lambda _{0}-\mu ^{-1},T\right) \left( \lambda _{0}I-T\right) -\mu
^{-1}R\left( \lambda _{0}-\mu ^{-1},T\right) =
\end{equation*}%
\begin{equation*}
=R\left( \lambda _{0}-\mu ^{-1},T\right) ((\lambda _{0}-\mu ^{-1})I-T+\mu
^{-1}I)-\mu ^{-1}R\left( \lambda _{0}-\mu ^{-1},T\right) =
\end{equation*}%
\begin{equation}{\label{equation:lb4}}
=I+\mu ^{-1}R\left( \lambda _{0}-\mu ^{-1},T\right) -\mu ^{-1}R\left( \lambda
_{0}-\mu ^{-1},T\right) =I
\end{equation}

Since, $R\left( \lambda _{0}-\mu ^{-1},T\right) $ comute with $T$ results
that%
\begin{equation*}
R\left( \lambda _{0}-\mu ^{-1},T\right) \left( \lambda _{0}I-T\right)
=\left( \lambda _{0}I-T\right) R\left( \lambda _{0}-\mu ^{-1},T\right) ,
\end{equation*}
which implies that%
\begin{equation*}
R(\lambda _{0},T)R\left( \lambda _{0}-\mu ^{-1},T\right) =R\left( \lambda
_{0}-\mu ^{-1},T\right) R(\lambda _{0},T).
\end{equation*}

From this equality and relation (\ref{equation:lb4}) results (\ref%
{equation:lb3}).

Moreover, 
\begin{equation*}
R\left( \lambda _{0}-\mu ^{-1},T\right) ,\left( \lambda _{0}I-T\right) \in
Q_{\mathcal{P}}(X),
\end{equation*}
so $R(\mu ,R(\lambda _{0},T))\in Q_{\mathcal{P}}(X)$, for all $\mid \mu \mid
>\epsilon ^{-1}$.

Therefore, $\sigma (Q_{\mathcal{P}},R(\lambda _{0},T))\subset D\left(
0,\epsilon ^{-1}\right).$

But if $\sigma (Q_{\mathcal{P}},R(\lambda _{0},T))$ is bounded then by
proposition \ref{proposition:gi2} there exists a family of seminorms $%
\mathcal{P}^{\prime }\in \mathcal{C}(X)$ which satisfies the condition

\begin{enumerate}
\item the families $\mathcal{P}=\left( p_{\alpha }\right) _{\alpha \in
\Lambda }$ and $\mathcal{P}^{\prime }=\left( p_{\alpha }^{\prime }\right)
_{\alpha \in \Lambda }$ has the same family of indices.

\item for every $\alpha \in \Lambda $ there exists $m_{\alpha },M_{\alpha
}>0 $ such that 
\begin{equation*}
m_{\alpha }p_{\alpha }\left( x\right) \leq p_{\alpha }^{\prime }\left(
x\right) \leq M_{\alpha }p_{\alpha }\left( x\right) ,\left( \forall \right)
x\in X;
\end{equation*}

\item $R(\lambda _{0},T)\in B_{\mathcal{P}}^{\prime }(X)$
\end{enumerate}

Moreover, it is obvious that $Q_{\mathcal{P}}^{\prime }(X)=Q_{\mathcal{P}%
}(X) $.

Let $\epsilon _{0}=\min \{\epsilon ,2^{-1}\left\Vert R(\lambda
_{0},T)\right\Vert _{\mathcal{P}^{\prime }}^{-1}\}$. Then for every $\mid
\mu \mid <\epsilon _{0}$ we have $\left\Vert \mu R(\lambda
_{0},T)\right\Vert _{\mathcal{P}^{\prime }}<\frac{1}{2}$, so the series $%
\sum\limits_{k=0}^{\infty }\left( -\mu \right) ^{k}R(\lambda _{0},T)^{k}$
converges in the algebra $B_{\mathcal{P}^{\prime }}(X)$ with respect to $%
\left\Vert \;\right\Vert _{\mathcal{P}^{\prime }}$. From the equalities%
\begin{equation*}
\left( I+\mu R(\lambda _{0},T)\right) \left( \sum\limits_{k=0}^{\infty
}\left( -\mu \right) ^{k}R(\lambda _{0},T)^{k}\right) =
\end{equation*}%
\begin{equation*}
=\sum\limits_{k=0}^{\infty }\left( -\mu \right) ^{k}R(\lambda
_{0},T)^{k}-\sum\limits_{k=0}^{\infty }\left( -\mu \right) ^{k+1}R(\lambda
_{0},T)^{k+1}=I
\end{equation*}
results that 
\begin{equation*}
\left( I+\mu R(\lambda _{0},T)\right) ^{-1}=\sum\limits_{k=0}^{\infty
}\left( -\mu \right) ^{k}R(\lambda _{0},T)^{k}
\end{equation*}

Using this relation and the following equalities 
\begin{equation*}
\sum\limits_{k=0}^{\infty }\left( -\mu \right) ^{k}R(\lambda
_{0},T)^{k+1}=\left( \sum\limits_{k=0}^{\infty }\left( -\mu \right)
^{k}R(\lambda _{0},T)^{k}\right) R(\lambda _{0},T)=
\end{equation*}%
\begin{equation*}
=\left[ I+\mu R(\lambda _{0},T)\right] ^{-1}R(\lambda _{0},T)=\left[
(\lambda _{0},T)\left( I+\mu R(\lambda _{0},T)\right) \right] ^{-1}=
\end{equation*}%
\begin{equation*}
=(\lambda _{0}I-T-\mu I)^{-1}=R(\mu +\lambda _{0},T)
\end{equation*}%
results that $\sum\limits_{k=0}^{\infty }\left( -\mu \right) ^{k}R(\lambda
_{0},T)^{k+1}$ converges in $B_{\mathcal{P}^{\prime }}(X)$ to $R(\mu
+\lambda _{0},T)$.

Since $\left\Vert \mu R(\lambda _{0},T)\right\Vert _{\mathcal{P}^{\prime }}<%
\frac{1}{2}$ results that 
\begin{equation*}
\left\Vert R(\mu +\lambda _{0},T)\right\Vert _{\mathcal{P}^{\prime }}\leq
2\left\Vert R(\lambda _{0},T)\right\Vert _{\mathcal{P}^{\prime }},\left(
\forall \right) \mid \mu \mid <\epsilon _{0}.
\end{equation*}%
so using lemma \ref{lemma:1} we have 
\begin{equation*}
p_{\alpha }^{\prime }(R(\mu +\lambda _{0},T)x)\leq 2\left\Vert R(\lambda
_{0},T)\right\Vert _{\mathcal{P}^{\prime }}p_{\alpha }^{\prime }\left(
x\right) ,
\end{equation*}%
for every $x\in X$ and every $\alpha \in $ $\Lambda $. The property (2) of
family $\mathcal{P}^{\prime }\in \mathcal{C}(X)$ implies that 
\begin{equation*}
p_{\alpha }(R(\mu +\lambda _{0},T)x)\leq m_{\alpha }^{-1}p_{\alpha }^{\prime
}(R(\mu +\lambda _{0},T)x)\leq
\end{equation*}%
\begin{equation*}
\leq m_{\alpha }^{-1}2\left\Vert R(\lambda _{0},T)\right\Vert _{\mathcal{P}%
^{\prime }}p_{\alpha }^{\prime }\left( x\right) \leq
\end{equation*}%
\begin{equation*}
\leq 2m_{\alpha }^{-1}M_{\alpha }\left\Vert R(\lambda _{0},T)\right\Vert _{%
\mathcal{P}^{\prime }}p_{\alpha }\left( x\right)
\end{equation*}%
for every $x\in X$ and every $\mid \mu \mid <\epsilon _{0}$. Therefore, 
\begin{equation*}
\hat{p}_{\alpha }(R(\mu +\lambda _{0},T)x)\leq 2m_{\alpha }^{-1}M_{\alpha
}\left\Vert R(\lambda _{0},T)\right\Vert _{\mathcal{P}^{\prime }}
\end{equation*}%
for every $\alpha \in $ $\Lambda $ and every $\mid \mu \mid <\epsilon _{0}$,
so the set $\{R(\lambda ,T){|}\lambda \in D\left( \lambda _{0},\epsilon
_{0}\right) \}$ is bounded in $Q_{\mathcal{P}}(X)$ and $\lambda _{0}\in \rho
_{W}(Q_{\mathcal{P}},T)$.
\end{proof}

\begin{proposition}
{\label{proposition:lb1}}If $(X,\mathcal{P})$ is a locally convex space such
that $\mathcal{P}\in \mathcal{C}_{0}(X)$, then for each operator $T\in 
\mathcal{LB}(X)$ there exists a family of seminorms $\mathcal{P}^{\prime
}\in \mathcal{C}(X)$ such that $\mathcal{P\approx P}^{\prime }$ and $T\in B_{%
\mathcal{P}^{\prime }}(X)$.
\end{proposition}
\begin{proof}
Since $T\in \mathcal{LB}(X)$ there exist $p_{0}\in \mathcal{P}$ such that  $m_{p_{0}q}(T)<\infty $, for each $q\in \mathcal{P}$ . From lemma \ref%
{lemma:op.semin.} results that 
\begin{equation}{\label{equation:lb5}}
q(Tx)\leq m_{p_{0}q}(T)p_{0}(x),~\left( \forall
\right) x\in X,\left( \forall \right) q\in \mathcal{P}
\end{equation}

For every $q\in \mathcal{P}$, we consider the application $q^{\prime
}:X\rightarrow R$ given by relation%
\begin{equation*}
q^{\prime }(x)=\max \left\{ q(x),m_{p_{0}q}(T)p_{0}(x)\right\} ,~\left(
\forall \right) x\in X.
\end{equation*}

It is easily to see that $\mathcal{P}^{\prime }=\{q^{\prime }|~q\in \mathcal{%
P}\}\in \mathcal{C}(X)$ and for every $q\in \mathcal{P}$ we have $q\leq
q^{\prime }$. Let $q^{\prime }\in \mathcal{P}^{\prime }$, where 
\begin{equation*}
q^{\prime }(x)=\max \left\{ q(x),m_{p_{0}q}(T)p_{0}(x)\right\} ,~\left(
\forall \right) x\in X.
\end{equation*}

Since $\mathcal{P}\in \mathcal{C}_{0}(X)$ results that there exists $%
q_{1}\in \mathcal{P}$ such that $p_{0}\leq q_{1}$ and $q\leq q_{1}$, so 
\begin{equation*}
q^{\prime }(x)=\max \left\{ q(x),m_{p_{0}q}(T)p_{0}(x)\right\} \leq
c_{q}q_{1}(x),~\left( \forall \right) x\in X,
\end{equation*}
where $c_{q}=\max \{1,m_{p_{0}q}(T)\}$. Therefore, for every $q^{\prime }\in 
\mathcal{P}^{\prime }$ there exists $q_{1}\in \mathcal{P}$ and $c_{q}>0$ such that $%
q^{\prime }\leq c_{q}q_{1}$. This will implies that $\mathcal{P\approx P}^{\prime
}$.

Moreover, from (\ref{equation:lb5}) results that for every $q^{\prime }\in 
\mathcal{P}^{\prime }$ we have 
\begin{equation*}
q^{\prime }(Tx)=\max \left\{ q(Tx),m_{p_{0}q}(T)p_{0}(Tx)\right\} \leq 
\end{equation*}%
\begin{equation*}
\leq \max
\left\{ m_{p_{0}q}(T)p_{0}(x),m_{p_{0}q}(T)\hat{p}_{0}(T)p_{0}(x)\right\} =
\end{equation*}%
\begin{equation*}
=m_{p_{0}q}(T)p_{0}(x)\max \left\{ 1,\hat{p}_{0}(T)\right\} \leq
c_{0}m_{p_{0}q}(T)q^{\prime }(x)
\end{equation*}
for all $x\in X$ (where $c_{0}=\max \left\{ 1,\hat{p}_{0}(T)\right\} $), so $T\in
B_{\mathcal{P}^{\prime }}(X)$.
\end{proof}

\begin{corollary}
{\label{corollary:lb1}}If $(X,\mathcal{P})$ is a locally convex space such
that $\mathcal{P}\in \mathcal{C}_{0}(X)$ and $\mathcal{LB}(X)=\mathcal{L}(X)$%
, then%
\begin{equation*}
\mathcal{LB}(X)=\mathcal{L}(X)=Q_{\mathcal{P}}(X)=(Q_{\mathcal{P}}(X))_{0}.
\end{equation*}
\end{corollary}
\begin{proof}
It is a direct consequence of proposition and propositions \ref%
{proposition:jo1}, \ref{proposition:jo2} and \ref{proposition:lb1}.
\end{proof}

\begin{corollary}
{\label{corollary:lb2}}If $(X,\mathcal{P})$ is a locally convex space such
that $\mathcal{P}\in \mathcal{C}_{0}(X)$ and $T\in \mathcal{LB}(X)$, then%
\begin{equation*}
\sigma (T)=\sigma _{lb}(T)=\sigma (Q_{\mathcal{P}},T)=\sigma (Q_{\mathcal{P}%
}^{0},T)=\sigma _{W}(Q_{\mathcal{P}},T).
\end{equation*}
\end{corollary}
\begin{proof}
The proposition \ref{proposition:spectru lb 1} give us the equality $\sigma
(T)=\sigma _{lb}(T)$. 

If $\lambda \in \rho _{lb}(T)$, then there exists a
scalar $\alpha $ and a locally bounded operator $S$ on $X$ such that $%
(\lambda I-T)^{-1}=\alpha I+S$. From proposition \ref{proposition:jo2} and %
\ref{proposition:lb1} results that $S\in (Q_{\mathcal{P}}(X))_{0}$, so $%
(\lambda I-T)^{-1}\in (Q_{\mathcal{P}}(X))_{0}$. Therefore, $\lambda \in
\rho (Q_{\mathcal{P}},T)$ and 
\begin{equation*}
\rho _{lb}(T)\subset \rho (Q_{\mathcal{P}%
}^{0},T)\subset \rho (Q_{\mathcal{P}},T)\subset \rho (T),
\end{equation*}
which prove that 
\begin{equation*}
\sigma (T)=\sigma _{lb}(T)=\sigma (Q_{\mathcal{P}},T)=\sigma (Q_{\mathcal{P}%
}^{0},T).
\end{equation*}

Moreover, from proposition \ref{proposition:spectru lb} results that the set 
$\sigma (Q_{\mathcal{P}},T)$ is closed. Proposition \ref{proposition:lb1}
implies that there exists some calibration $\mathcal{P}^{\prime }\in 
\mathcal{C}(X)$ such that $\mathcal{P\approx P}^{\prime }$ and $T\in B_{%
\mathcal{P}^{\prime }}(X)$. Since $Q_{\mathcal{P}}(X)=Q_{\mathcal{P}^{\prime
}}(X)$ (proposition \ref{proposition:jo1}) and $\sigma (Q_{\mathcal{P}%
},T)=\sigma (Q_{\mathcal{P}^{\prime }},T)$ is closed, the theorem \ref%
{theorem:lb1} will prove that 
\begin{center}
$\sigma (Q_{\mathcal{P}},T)=\sigma (Q_{\mathcal{P}^{\prime }},T)=\sigma _{W}(Q_{\mathcal{P}^{\prime}},T)=\sigma _{W}(Q_{\mathcal{P}},T)$.
\end {center}
\end{proof}

\begin{remark}
If $(X,\mathcal{P})$ is a locally convex space such that $\mathcal{P}\in 
\mathcal{C}_{0}(X)$, then the corollary \ref{corollary:lb2} show that for
locally bounded operators on X the functional calculus presented in the
previous section is a natural generalization of the functional calculus for
bounded operators on Banach spaces.
\end{remark}

\end{document}